\documentclass{article}
\usepackage{amsmath,amsthm,amssymb}

\newcommand{\deopN}{\de
\hspace{-.3ex}\raisebox{1.1ex}[0pt][0pt]{\scriptsize\fontshape{n}\selectfont
op}\hspace{-1.8ex}\raisebox{-0.2ex}[0pt][0pt]{$_N$}\hspace{1.8ex}}
\newcommand{\ltimesred}{\;\mbox{$_\text{\rm r}$}\hspace{-0.8ex}\ltimes}
\newcommand{\red}{_\text{\rm r}}
\newcommand{\K}{\mathcal{K}}
\newcommand{\M}{\operatorname{M}}

\newcommand{\Ah}{\hat{A}}
\newcommand{\Atw}{A_\tw}

\newcommand{\etao}{\eta^1}
\newcommand{\util}{\tilde{u}}
\newcommand{\vtil}{\tilde{v}}
\newcommand{\wtil}{\tilde{w}}
\newcommand{\pitil}{\tilde{\pi}}
\newcommand{\zetatil}{\tilde{\zeta}}
\newcommand{\psitil}{\tilde{\psi}}
\newcommand{\Gatil}{\tilde{\Gamma}}
\newcommand{\Stil}{\tilde{S}}
\newcommand{\sitil}{\tilde{\sigma}}
\newcommand{\Util}{\tilde{U}}
\newcommand{\rhotil}{\tilde{\rho}}

\newcommand{\T}{\mathbb{T}}

\newcommand{\nabp}{\nab \hspace{-1.05ex}
\rule[.5ex]{.2ex}{.8ex}   \hspace{1.05ex}}

\newcommand{\Ad}{\operatorname{Ad}}
\newcommand{\sitiltil}{\tilde{\si}\hspace{-.7ex}\tilde{\rule{0ex}{1.4ex}}\hspace{.7ex}}

\newcommand{\nabtil}{\tilde{\nabla}}
\newcommand{\tetiltil}{\tilde{\te}
\hspace{-.4ex}\tilde{\rule{0ex}{2.0ex}}\hspace{.4ex}}
\newcommand{\vfitiltil}{\tilde{\vfi}
\hspace{-.6ex}\tilde{\rule{0ex}{1.4ex}}\hspace{.6ex}}

\newcommand{\nabte}{\nabla_\te}
\newcommand{\na}{\circ}

\newcommand{\tauo}{\tau^1}
\newcommand{\taut}{\tau^2}
\newcommand{\B}{{\rm B}}

\newcommand{\tw}{{\operatorname{m}}}
\newcommand{\detw}{\de_\tw}
\newcommand{\Jtw}{J_\tw}
\newcommand{\Jhtw}{\Jh_\tw}
\newcommand{\sdetw}{\delta_\tw}
\newcommand{\Wtw}{W_\tw}
\newcommand{\Mtw}{M_\tw}
\newcommand{\Rtw}{R_\tw}
\newcommand{\tautw}{\tau^\tw}
\newcommand{\nutw}{\nu_\tw}
\newcommand{\sdehtw}{\hat{\delta}_\tw}
\newcommand{\vfitw}{\vfi_\tw}
\newcommand{\detwop}{\de
\hspace{-.3ex}\raisebox{0.9ex}[0pt][0pt]{\scriptsize\fontshape{n}\selectfont
op}\hspace{-1.6ex}_\tw\hspace{1.6ex}}
\newcommand{\latw}{\Lambda_\tw}
\newcommand{\cWtw}{\mathcal{W}_\tw}

\newcommand{\Vtw}{V_\tw}
\newcommand{\psitw}{\psi_\tw}
\newcommand{\Utw}{U_\tw}
\newcommand{\Ptw}{P_\tw}
\newcommand{\Mhtw}{\hat{M}_\tw}
\newcommand{\dehtw}{\hat{\Delta}_\tw}
\newcommand{\Ahtw}{\hat{A}_\tw}

\newcommand{\tauh}{\hat{\tau}}

\newcommand{\Xtil}{\tilde{X}}

\newcommand{\Mhh}{\hat{M}
\hspace{-1.05ex}\hat{\rule{0ex}{2.0ex}}\hspace{1.05ex}}
\newcommand{\dehh}{\hat{\de}
\hspace{-.95ex}\hat{\rule{0ex}{2.05ex}}\hspace{.95ex}}

\newcommand{\cC}{\mathcal{C}}

\newcommand{\Wtil}{\tilde{W}}

\newcommand{\kruisje}[1]{\, \mbox{$_{#1}$}\hspace{-.2ex}\mbox{$\ltimes$} \,}

\newcommand{\Vtil}{\tilde{V}}
\newcommand{\alh}{\hat{\alpha}}
\newcommand{\Mh}{\hat{M}}
\newcommand{\Si}{\Sigma}

\newcommand{\dpr}{^{\prime\prime}}

\newcommand{\late}{\Lambda_\theta}
\newcommand{\te}{\theta}
\newcommand{\tetil}{\tilde{\theta}}

\newcommand{\Nte}{\cN_\theta}
\newcommand{\Tr}{{\operatorname{Tr}}}
\newcommand{\nsf}{n.s.f.\ }
\newcommand{\Jtil}{\tilde{J}}
\newcommand{\Jte}{J_\te}
\newcommand{\alhh}{\hat{\alpha}\hspace{-.7ex}\hat{\rule{0ex}{1.45ex}}\hspace{.7ex}}
\newcommand{\alhhklein}{\hat{\alpha}\hspace{-.55ex}\hat{\rule[-.4ex]{0ex}{1.45ex}}\hspace{.55ex}}
\newcommand{\Mo}{M_1}
\newcommand{\Mt}{M_2}
\newcommand{\deo}{\de_1}
\renewcommand{\det}{\de_2}
\newcommand{\deht}{\hat{\de}_2}
\newcommand{\Mht}{\hat{M}_2}

\newcommand{\Wh}{\hat{W}}

\newcommand{\Nfih}{\cN_{\vfih}}

\newcommand{\Nh}{\hat{N}}

\newcommand{\Rh}{\hat{R}}

\newcommand{\R}{\mathbb{R}}

\newcommand{\C}{\mathbb{C}}

\newcommand{\Ga}{\Gamma}
\newcommand{\deh}{\hat{\Delta}}

\newcommand{\vfih}{\hat{\vfi}}

\newcommand{\lah}{\hat{\la}}

\newcommand{\pih}{\hat{\pi}}
\newcommand{\sih}{\hat{\si}}
\newcommand{\sdeh}{\hat{\delta}}
\newcommand{\cU}{{\cal U}}

\newcommand{\nab}{\nabla}
\newcommand{\Jh}{\hat{J}}
\newcommand{\nabh}{\hat{\nab}}
\newcommand{\cZ}{{\cal Z}}

\newcommand{\cN}{{\cal N}}
\newcommand{\cM}{{\cal M}}

\newcommand{\ot}{\otimes}
\newcommand{\sla}{\lambda}
\newcommand{\la}{\Lambda}

\newcommand{\om}{\omega}
\newcommand{\io}{\iota}
\newcommand{\vfi}{\varphi}

\newcommand{\al}{\alpha}
\newcommand{\be}{\beta}

\newcommand{\sde}{\delta}
\newcommand{\de}{\Delta}

\newcommand{\si}{\sigma}
\newcommand{\Mfi}{{\cal M}_{\vfi}}
\newcommand{\Nfi}{{\cal N}_{\vfi}}
\newcommand{\Mpsi}{{\cal M}_{\psi}}
\newcommand{\Npsi}{{\cal N}_{\psi}}

\newcommand{\cst}{\text{C}$\hspace{0.1mm}^*$}

\newcommand{\Goop}{G_1
\hspace{-1ex}\raisebox{1ex}[0pt][0pt]{\scriptsize\fontshape{n}\selectfont
op}}
\newcommand{\deop}{\de \hspace{-.3ex}\raisebox{0.9ex}[0pt][0pt]{\scriptsize\fontshape{n}\selectfont op}}
\newcommand{\deopo}{\de_1
\hspace{-1ex}\raisebox{1ex}[0pt][0pt]{\scriptsize\fontshape{n}\selectfont
op}}
\newcommand{\deoop}{\deopo}
\newcommand{\deopt}{\de_2 \hspace{-1ex}\raisebox{1ex}[0pt][0pt]{\scriptsize\fontshape{n}\selectfont op}}
\newcommand{\detop}{\deopt}

\newcommand{\dehop}{\deh \hspace{-.3ex}\raisebox{0.9ex}[0pt][0pt]{\scriptsize\fontshape{n}\selectfont op}}
\newcommand{\dehopo}{\deh_1
\hspace{-1ex}\raisebox{1ex}[0pt][0pt]{\scriptsize\fontshape{n}\selectfont
op}}
\newcommand{\dehoop}{\dehopo}
\newcommand{\dehopt}{\deh_2 \hspace{-1ex}\raisebox{1ex}[0pt][0pt]{\scriptsize\fontshape{n}\selectfont op}}

\newcommand{\Op}{\raisebox{0.9ex}[0pt][0pt]{\scriptsize\fontshape{n}\selectfont op}\,}

\newcommand{\recht}{\rightarrow}

\numberwithin{equation}{section}

{\theoremstyle{definition}\newtheorem{definitiona}{Definition}[section]
\newtheorem{notation}[definitiona]{Notation}

}
\newtheorem{proposition}[definitiona]{Proposition}
\newtheorem{lemma}[definitiona]{Lemma}
\newtheorem{theorem}[definitiona]{Theorem}

\newenvironment{definition}{\begin{definitiona}}{\mbox{} \hfill
    $\blacktriangle$ \end{definitiona}}

\newcommand{\eindeformule}{\mbox{}\hfill\raisebox{0.6cm}[0mm][0mm]{$\blacktriangle$}\vspace{-0.6cm}}

\begin{document}
\begin{center}
{\Large\bf Double crossed products of locally compact quantum groups}

\bigskip

{\sc by Saad Baaj and Stefaan
  Vaes}
\end{center}

{\footnotesize Laboratoire de Math{\'e}matiques Pures;
Universit{\'e} Blaise Pascal; B{\^a}timent de Math{\'e}matiques; F--63177 Aubi{\`e}re
Cedex (France) \\ e-mail: Saad.Baaj@math.univ-bpclermont.fr \vspace{0.3ex}\\
Institut de Math{\'e}matiques de Jussieu; Alg{\`e}bres d'Op{\'e}rateurs et
Repr{\'e}sentations; 175, rue du Chevaleret; F--75013 Paris (France) \\
e-mail: vaes@math.jussieu.fr}

\bigskip

\begin{abstract}\noindent
For a matched pair of locally compact quantum groups, we construct the double crossed product as a locally compact quantum group. This construction
generalizes Drinfeld's quantum double construction. We study \cst-algebraic properties of these double crossed products and several links between
double crossed products and bicrossed products. In an appendix, we study the Radon-Nikodym derivative of a weight under a quantum group action,
following Yamanouchi and obtain, as a corollary, a new characterization of closed quantum subgroups.
\end{abstract}

\section{Introduction}
In this paper, we study the so-called \emph{double crossed product}
construction in the framework of \emph{locally compact quantum groups},
generalizing Drinfeld's quantum double construction and bringing it
into the topological setting.

The history of locally compact quantum groups goes back to the 1960's, when G. Kac tried to answer the following basic question: is it possible to
unify locally compact (non-abelian) groups and their dual objects into one theory? A positive answer was given by Kac, Kac \& Vainerman and
independently, by Enock \& Schwartz. They defined what Enock \& Schwartz called \emph{Kac algebras}, see \cite{E-S} for an overview. In the 1980's,
several new examples of objects certainly to be considered as quantum groups, were constructed and were shown not to be Kac algebras. Henceforth,
many efforts were made to enlarge the definition of Kac algebras to include these new examples and in particular, to include Woronowicz'
\emph{compact quantum groups}, see \cite{Wor1,Wor2}. Many important contributions have been made by Skandalis \& the first author \cite{BS},
Woronowicz \cite{Wor3}, Masuda \& Nakagami \cite{Mas-Nak} and Van Daele \cite{VD}. Finally, the theory of locally compact quantum groups was proposed
by Kustermans and the second author in \cite{KV1,KV2}.

One of the most celebrated constructions in the early, algebraic theory of quantum groups is \emph{Drinfeld's quantum double construction}
\cite{Drin}. More generally, Majid \cite{Maj} defined the double crossed product for mutually matched Hopf algebras, yielding Drinfeld's quantum
double when a Hopf algebra is matched with its dual in the canonical way.

It is a natural idea to perform analogous constructions in the
topological framework of locally compact (l.c.) quantum groups. In
fact, the double crossed product construction was developed for
multiplicative unitaries in \cite{BS} (it is called \emph{$Z$-produit
  tensoriel}, there). Nevertheless, some extra conditions had to be
imposed in order to obtain the double crossed product and it was not
clear if they were satisfied for all matchings of two multiplicative
unitaries.

In this paper, we show that any matching of l.c.\ quantum groups
allows to perform the double crossed product construction and to
obtain again a l.c.\ quantum group. In particular, the quantum double
of a l.c.\ quantum group is again a l.c.\ quantum group: a result
already shown for Woronowicz algebras by Yamanouchi \cite{Yam3}.

The notion of a matching of l.c.\ quantum groups was developed by
Vainerman and the second author in \cite{VV} and the \emph{bicrossed
  product} was constructed. We explain all these different
constructions in a more detailed way for a matching of ordinary l.c.\ groups. One says that a pair of l.c.\ groups $G_1,G_2$ is a matched pair, if
there exists a large group $G$, such that $G_1,G_2$ are closed subgroups, satisfying $G_1 \cap G_2 = \{e\}$ and such that the complement of $G_1 G_2$
has Haar measure zero in $G$. In fact, this allows to define an action $(\al_g)_{g \in G_1}$ of $G_1$ by (almost everywhere defined) transformations
of the measure space $G_2$ (not by automorphisms) and an action $(\be_s)_{s \in G_2}$ of $G_2$ by transformations of the measure space $G_1$ such
that $\al_g(s) \; \be_s(g) = gs$ almost everywhere. Conversely, one can say that $G_1,G_2$ are matched when they act on each other in such a
compatible way. Such a compatible pair of actions allows to reconstruct the l.c.\ group $G$, this is the \emph{double crossed product}, and to
construct a l.c.\ quantum group with underlying von Neumann algebra $G_1 \kruisje{\al} L^\infty(G_2)$, this is the \emph{bicrossed product}.

In the quantum setting, a matching of two l.c.\ quantum groups was
defined in \cite{VV} as a compatible pair of actions of the two quantum
groups on each other. So, it is a natural question if one can find the
double crossed product, i.e.\ a big l.c.\ quantum group such that the
matched quantum groups are closed quantum subgroups and such that the matching
is determined by this inclusion in the double crossed product, as
above by the equation $\al_g(s) \; \be_s(g) = gs$. In our main result,
Theorem \ref{thm.double} below, we show that there always exists a double
crossed product quantum group and we calculate its structural
ingredients. In Proposition \ref{prop.characterization} below, we show
that the notions of a \emph{matched pair} and a \emph{quantum group
  with two closed quantum subgroups} are in bijective correspondence.

In \cite{BSV}, Skandalis and the authors studied the \cst-algebraic features of the bicrossed product associated with a matched pair $G_1,G_2 \subset
G$ of l.c.\ groups. The main result was that the bicrossed product is \emph{regular} if and only if $G = G_1G_2$ homeomorphically, while it is
\emph{semi-regular} if and only if $G_1G_2$ is an open subset of $G$. Examples of non-semi-regular l.c.\ quantum groups were given using examples of
matched pairs $G_1, G_2 \subset G$ such that $G_1G_2$ has empty interior in $G$ (but, as needed, $G_1G_2$ is of full Haar measure in $G$). The fact
that $G=G_1 G_2$ homeomorphically can be expressed in \cst-algebraic terms as $C_0(G) = C_0(G_1) \ot C_0(G_2)$, while the fact that $G_1 G_2$ is an
open subset of $G$ becomes $C_0(G_1) \ot C_0(G_2) \subset C_0(G)$. In Section \ref{sec.cstar}, we study the same kind of \cst-algebraic questions in
the quantum setting and provide, in particular, the same kind of necessary and sufficient condition for the bicrossed product to be (semi-)regular.
More specifically, regularity is equivalent with the equation $\Atw = A_1 \ot A_2$ and semi-regularity with $A_1 \ot A_2 \subset \Atw$, where
$A_1,A_2$ and $\Atw$ denote the \cst-algebras associated with the two matched quantum groups and the double crossed product, respectively.

The structure of the paper is as follows: we first recall the necessary preliminaries in Section \ref{sec.prelim}, in particular the basics of the
theory of l.c.\ quantum groups and the main results on matchings and bicrossed products from \cite{VV}. Sections \ref{sec.doubleone} --
\ref{sec.doublelc} are the core of the paper: we define the double crossed product for a matched pair and prove that it is a l.c.\ quantum group. In
order to do so, we have to produce the Haar weights and we will need a continuous interplay between the double crossed product and the bicrossed
products, about which we already know a lot from \cite{VV}. In the short Section \ref{sec.actions}, we show how actions of double crossed products on
von Neumann algebras and actions of bicrossed products are related by an intermediate von Neumann algebra in the Jones tower. In Section
\ref{sec.char}, we prove the bijective correspondence between a matched pair and a quantum group with two closed quantum subgroups in a good position
w.r.t.\ each other. The example of a (generalized) quantum double construction is treated in Section \ref{sec.doublegen} and we deal with the
\cst-algebraic features of our constructions in Section \ref{sec.cstar}.

One of our major tools is the \emph{Radon-Nikodym derivative} of a weight under a l.c.\ quantum group action, as developed originally by Yamanouchi
\cite{Yam2}. In an appendix, Section \ref{sec.app}, we give a simplified account of Yamanouchi's interesting results. As an application, somehow
aside of the main line of our paper, we give a characterization of \emph{closed quantum subgroups} as von Neumann subalgebras invariant under the
comultiplication and the antipode.

{\bf Acknowledgment } The authors would like to thank Georges Skandalis for many fruitful discussions on the topic of this paper.

\section{Preliminaries} \label{sec.prelim}

We use the following conventions: $\ot$ will denote, depending on the case, the tensor product of Hilbert spaces, von Neumann algebras or
\cst-algebras. There should be no confusion. The domain of an unbounded mapping $S$ is denoted by $D(S)$. We write $\si$ for the flip map on the
tensor product $M \ot M$ of two von Neumann algebras and $\Si$ for the flip map on the tensor product of two Hilbert spaces $H \ot H$.

\subsection{A little bit of weight theory}

Let $\vfi$ be a normal, semi-finite, faithful (n.s.f.) weight on a von Neumann algebra $M$. Then, we write
$$\cM_\vfi^+ = \{ x \in M^+ \mid \vfi(x) < \infty \} \quad\text{and} \quad \cN_\vfi = \{ x \in M \mid \vfi(x^* x) < \infty \} \; .$$
We represent $M$ on its GNS-space $H$ and obtain a GNS-map $\la : \Nfi \recht H$. Denote by $(\si_t)$ the modular group of $\vfi$.

Suppose that $k$ is a strictly positive, self-adjoint operator affiliated with $M$ and satisfying $\si_t(k) = \sla^{t} k$ for all $t \in \R$ and a
number $\sla > 0$. Then, we can define a new n.s.f.\ weight $\vfi_k$ which is formally given by $\vfi_k(x) = \vfi(k^{1/2} x k^{1/2})$. For a precise
definition, we refer to \cite{V2}. We also get a canonical GNS-map $\la_{\vfi_k} : \cN_{\vfi_k} \recht H$ for the weight $\vfi_k$, which is formally
given by $\la_{\vfi_k}(x) = \la(x k^{1/2})$.

\subsection{Locally compact quantum groups}
Our references for the theory of l.c.\ quantum groups are \cite{KV1,KV2}.
\begin{definitiona}
A pair $(M,\de)$ is called a (von Neumann algebraic) l.c.\ quantum group when
\begin{itemize}
\item $M$ is a von Neumann algebra and $\de : M \recht M \ot M$ is
a normal and unital $*$-homomorphism satisfying the coassociativity relation : $(\de \ot \io)\de = (\io \ot \de)\de$.
\item There exist n.s.f.\ weights $\varphi$ and $\psi$ on $M$ such that
\begin{itemize}
\item $\varphi$ is left invariant in the sense that $\varphi \bigl( (\om \ot
\io)\de(x) \bigr) = \varphi(x) \om(1)$ for all $x \in \Mfi^+$ and $\om \in M_*^+$,
\item $\psi$ is right invariant in the sense that $\psi \bigl( (\io \ot
\om)\de(x) \bigr) = \psi(x) \om(1)$ for all $x \in \Mpsi^+$ and $\om \in M_*^+$.
\end{itemize}
\end{itemize} \eindeformule
\end{definitiona}
Fix such a l.c.\ quantum group $(M,\de)$.

Represent $M$ in the GNS-construction of $\vfi$ with GNS-map $\la : \Nfi \recht H$. We define a unitary $W$ on $H \ot H$ by
$$W^* (\Lambda(a) \ot \Lambda(b)) = (\Lambda \ot \Lambda)(\de(b)(a \ot 1)) \quad\text{for all}\; a,b \in N_{\phi}
\; .$$ Here, $\Lambda \ot \Lambda$ denotes the canonical GNS-map for the tensor product weight $\varphi \ot \varphi$. One proves that $W$ satisfies
the pentagonal equation: $W_{12} W_{13} W_{23} = W_{23} W_{12}$, and we say that $W$ is a \emph{multiplicative unitary}. It is the \emph{left regular
representation}. The von Neumann algebra $M$ is the strong closure of the algebra $\{ (\io \ot \om)(W) \mid \om \in \B(H)_* \}$ and $\de(x) = W^* (1
\ot x) W$, for all $x \in M$. Next, the l.c.\ quantum group $(M,\de)$ has an antipode $S$, which is the unique $\si$-strong$^*$ closed linear map
from $M$ to $M$ satisfying $(\io \ot \om)(W) \in D(S)$ for all $\om \in B(H)_*$ and $S(\io \ot \om)(W) = (\io \ot \om)(W^*)$ and such that the
elements $(\io \ot \om)(W)$ form a $\si$-strong$^*$ core for $S$. The antipode $S$ has a polar decomposition $S = R \tau_{-i/2}$ where $R$ is an
anti-automorphism of $M$ and $(\tau_t)$ is a strongly continuous one-parameter group of automorphisms of $M$. We call $R$ the \emph{unitary antipode}
and $(\tau_t)$ the \emph{scaling group} of $(M,\de)$. From \cite{KV1}, Proposition~5.26 we know that $\sigma (R \ot R) \de = \de R$. So $\varphi R$
is a right invariant weight on $(M,\de)$ and we take $\psi:= \varphi R$.

Let us denote by $(\sigma_t)$ the modular automorphism group of $\varphi$. From \cite{KV1}, Proposition~6.8 we get the existence of a number $\nu >
0$, called the \emph{scaling constant}, such that $\psi \, \sigma_t = \nu^{-t} \, \psi$ for all $t \in \R$. Hence, we get the existence of a unique
positive, self-adjoint operator $\sde$ affiliated to $M$, such that $\sigma_t(\sde) = \nu^t \, \sde$ for all $t \in \R$ and $\psi = \varphi_{\sde}$,
see \cite{KV1}, Definition~7.1. The operator $\sde$ is called the \emph{modular element} of $(M,\de)$ and the canonically associated GNS-map is
denoted by $\Ga : \Npsi \recht H$. Because of the right invariance of $\psi$, we find another multiplicative unitary, denoted by $V$ and called the
\emph{right regular representation}:
$$V (\Gamma(x) \ot \Gamma(y)) = (\Gamma \ot \Gamma)(\de(x) (1 \ot y)) \quad\text{for all}\quad x,y \in \Npsi \; .$$
The scaling constant can be characterized as well by the relative invariance $\varphi \, \tau_t = \nu^{-t} \, \varphi$.

The \emph{dual l.c.\ quantum group} $(\Mh,\deh)$ is defined in \cite{KV1}, Section~8. Its von Neumann algebra $\Mh$ is the strong closure of the
algebra $\{(\om \ot \io)(W) \mid \om \in \B(H)_* \}$ and the comultiplication is given by $\deh(x) = \Sigma W (x \ot 1) W^* \Sigma$ for all $x \in
\Mh$. On $\Mh$ there exists a canonical left invariant weight $\vfih$ with GNS-map $\lah : \Nfih \recht H$ and the associated multiplicative unitary
is denoted by $\hat{W}$. From \cite{KV1}, Proposition~8.16, it follows that $\hat{W} = \Sigma W^* \Sigma$.

Since $(\Mh,\deh)$ is again a l.c.\ quantum group, we can introduce the antipode $\hat{S}$, the unitary antipode $\hat{R}$, the scaling group
$(\hat{\tau}_t)$ and the modular element $\sdeh$ exactly as we did it for $(M,\de)$. Also, we can again construct the dual of $(\Mh,\deh)$, starting
from the left invariant weight $\vfih$ with GNS-construction $(H,\io,\lah)$. From \cite{KV1}, Theorem~8.29 we have that the bidual l.c.\ quantum
group $(\Mhh,\dehh)$ is isomorphic to $(M,\de)$.

We denote by $(\sih_t)$ the modular automorphism groups of the weight $\vfih$. The modular conjugations of the weights $\varphi$ and $\hat\varphi$
will be denoted by $J$ and $\Jh$ respectively. Then it is worthwhile to mention that $$R(x) = \Jh x^* \Jh \quad\text{for all} \; x \in M
\qquad\text{and}\qquad \Rh(y) = J y^* J \quad\text{for all}\; y \in \Mh \; .$$ Some important relations involving $J$ and $\Jh$ are
$$
\Jh J = \nu^{i/4} J \Jh \; , \quad W^* = (\Jh \ot J) W (\Jh \ot J) \; , \quad V = (\Jh \ot \Jh)\Si W^* \Si (\Jh \ot \Jh) \; .
$$
The modular operators of $\vfi$ and $\vfih$ are denoted by $\nab$ and $\nabh$, respectively. Finally, we introduce the implementation $P$ of
$(\tau_t)$ as the unique strictly positive, self-adjoint operator on $H$ such that $P^{it} \la(x) = \nu^{t/2} \la(\tau_t(x))$. We mention that
$P=\hat{P}$, i.e., $P$ is as well the implementation of $(\tauh_t)$. For all kinds of relations between the operators on $H$ introduced so far, we
refer to Proposition 2.4 in \cite{VVD}.

Every l.c.\ quantum group has an \emph{associated \cst-algebra} of \lq continuous functions tending to zero at infinity\rq\ and we denote them by $A$
and $\Ah$:
$$A := [(\io \ot \om)(W) \mid \om \in \B(H)_* ] \; , \quad \Ah = [(\om
\ot \io)(W) \mid \om \in \B(H)_*] \; ,$$
where we use the following notation.

\begin{notation}
When $X$ is a subset of a Banach space, we denote by $[X]$ the closed
linear span of $X$. We sometimes use a notation like $[(\io \ot
\om)(W)]$ instead of $[(\io \ot \om)(W) \mid \om \in \B(H)_* ]$,
because $\om$ will always run through $\B(H)_*$ for the appropriate
Hilbert space $H$.
\end{notation}

For any multiplicative unitary $W$ on a Hilbert space $H$, we introduce (following
\cite{B,BS}) the algebra $\cC(W)$ by the formula
$$\cC(W) := \{ (\io \ot \om)(\Si W) \mid \om \in \B(H)_* \} \; .$$
We say that $(M,\de)$ is \emph{regular} (see \cite{BS}) if $[\cC(W)] = \K(H)$, the compact operators on $H$. We say that $(M,\de)$ is
\emph{semi-regular} (see \cite{B}) if $\K(H) \subset [\cC(W)]$. From Proposition 2.6 in \cite{BSV}, we get that
$$[\cC(W)] = [JAJ \; \Ah] \; .$$
Not all l.c.\ quantum groups are semi-regular, see \cite{BSV}.

\subsection{Actions of l.c.\ quantum groups}

We use \cite{V} as a reference for actions of l.c.\ quantum groups, but we repeat the necessary elements of the theory. Let $(M,\de)$ be a l.c.\
quantum group and $N$ a von Neumann algebra. A faithful, normal $^*$-homomorphism $\al : N \recht M \ot N$ is called a {\it (left) action} of
$(M,\de)$ on $N$ if $(\de \ot \io)\al = (\io \ot \al)\al$. One can define the crossed product as
$$M \kruisje{\al} N = (\al(N) \cup \Mh \ot 1)\dpr \in \B(H) \ot N \; .$$
Given a n.s.f.\ weight $\te$ on $N$, one can define the \emph{dual weight} $\tetil$ on the crossed product $M \kruisje{\al} N$. Representing $N$ on
the GNS-space $H_\te$ of $\te$, we get a GNS-map $\la_\te : \Nte \recht H_\te$. From \cite{V}, Definition~3.4 and Proposition~3.10, we get a
canonical GNS-construction for $\tetil$ on the Hilbert space $H \ot H_\te$, with GNS-map determined by
$$(x \ot 1) \al(y) \mapsto \lah(x) \ot \la_\te(y) \quad\text{for}\quad
x \in \Nfih, y \in \Nte$$ and where these elements $(x \ot 1) \al(y)$ span a core for the GNS-map, when we equip $M \kruisje{\al} N$ with its
$\si$-strong$^*$ topology and $H \ot H_\te$ with the norm topology.

If we denote by $J_\te$ the modular conjugation of $\te$ and by $\Jtil_\te$ the modular conjugation of the dual weight $\tetil$ in the canonical
GNS-construction for $\tetil$, we can define the unitary $U_\al = \Jtil_\te (\Jh \ot J_\te)$, which is called the {\it unitary implementation} of
$\al$, see \cite{V}. It satisfies
$$U_\al \in M \ot \B(H_\te) \; , \quad (\de \ot \io)(U_\al) = U_{\al,23} U_{\al,13} \quad\text{and}\quad \al(x) = U_\al (1 \ot x) U_\al^*
\;\;\text{for}\; x \in N \; .$$ Because $U_\al^*$ is a corepresentation of $(M,\de)$, we get that $U_\al \in \M(A \ot \K(H_\te))$, where $A \subset
M$ is the \cst-algebra associated with $(M,\de)$.

If $\al : N \recht M \ot N$ is an action, the \emph{fixed point algebra} is denoted by $N^\al$ and defined as the von Neumann subalgebra of elements
$x \in N$ satisfying $\al(x) = 1 \ot x$.

Up to now, we dealt with left actions, but this is only a matter of
convention. We can (and will) also consider right actions $\al : N
\recht N \ot M$ of $(M,\de)$ on a von Neumann algebra $N$, which means
that $(\al \ot \io)\al = (\io \ot \de)\al$.

\subsection{Bicrossed products of l.c.\ quantum groups}

We give a quick overview of some of the results of \cite{VV}. Fix two l.c.\ quantum groups $(M_1,\de_1)$ and $(M_2,\de_2)$.

\begin{definitiona} \label{defin.matching}
A normal, faithful $^*$-homomorphism $\tw : M_1 \ot M_2 \recht M_1 \ot M_2$ is called a {\it matching} of $(M_1,\de_1)$ and $(M_2,\de_2)$ if $$(\de_1
\ot \io)\tw = \tw_{23} \tw_{13} (\de_1 \ot \io) \quad\text{and}\quad
(\io \ot \de_2)\tw = \tw_{13}\tw_{12}(\io \ot \de_2) \; .$$ \eindeformule
\end{definitiona}

We fix a matching $\tw$ of $(M_1,\de_1)$ and $(M_2,\de_2)$. As in \cite{VV}, we define
$$\al : M_2 \recht M_1 \ot M_2 : \al(x) = \tw(1 \ot x) \quad\text{and}\quad \be : M_1 \recht M_1 \ot M_2 : \be(y) = \tw(y \ot 1) \; .$$
Then, $\al$ is a left action of $(M_1,\de_1)$ on the von Neumann algebra $M_2$, while $\be$ is a right action of $(M_2,\detop)$ on $M_1$. So, we can
define two crossed product von Neumann algebras $M$ and $\Mh$ on $H_1 \ot H_2$:
$$M = M_1 \kruisje{\al} M_2 = \bigl(\al(M_2) \cup \Mh_1 \ot 1)\dpr \quad\text{and}\quad \Mh = M_1 \rtimes_\be M_2 = \bigl( \be(M_1) \cup 1 \ot \Mh_2
\bigr)\dpr \; .$$ In \cite{VV}, it is proven that we can define a l.c.\ quantum group structure $(M,\de)$ on the von Neumann algebra $M$ such that
$\Mh$ is the dual quantum group. The comultiplications on $M$ and its dual $\Mh$ are given by
\begin{align*}
\de \al &= (\al \ot \al)\de_2 \; , \quad (\io \ot \de)(W_1 \ot 1) = W_{1,14} \; ((\io \ot \al)\be \ot \io)(W_1)_{1453} \; , \\
\deh \be &= (\be \ot \be)\de_1 \; , \quad (\deh \ot \io)(1 \ot \Wh_2) = (\io \ot (\be \ot \io)\al)(\Wh_2)_{2345} \; \Wh_{2,25} \; .
\end{align*}
The l.c.\ quantum group $(M,\de)$ is called the \emph{bicrossed product} associated with the matching $\tw$ of $(\Mo,\deo)$ and $(\Mt,\det)$.

Define $\vfi$ to be the dual weight on $M$ of the weight $\vfi_2$ on $M_2$, with canonical GNS-map determined by
$$\la\bigl( (a \ot 1)\al(x) \bigr) = \lah_1(a) \ot \la_2(x) \quad\text{for}\;\; a \in \cN_{\vfih_1} , x \in \cN_{\vfi_2} \; .$$
Then, $\vfi$ is the left invariant weight on $(M,\de)$ and the left regular representation is given by
$$W = (\al \ot \io)(W_2)_{124} \; (\io \ot \be)(\Wh_1)_{134} \; .$$
Observe that, on $H_1 \ot H_2$, we have the modular conjugations $J$ and $\Jh$ of the invariant weights $\vfi$ and $\vfih$. So, we also have
\begin{align*}
U_\al & = J(\Jh_1 \ot J_2) \in M_1 \ot \B(H_2) \; , \quad \al(x) = U_\al(1 \ot x)U_\al^* \; , \\ U_\be &= \Jh(J_1 \ot \Jh_2) \in \B(H_1) \ot M_2 \; ,
\quad  \be(y) = U_\be (y \ot 1) U_\be^* \; .
\end{align*}

\begin{notation} \label{nota.crucial}
Throughout Sections \ref{sec.doubleone} -- \ref{sec.actions}, we will fix a matching $\tw$ of two l.c.\ quantum groups $(M_1,\de_1)$ and
$(M_2,\de_2)$. All notations with indices $1,2$ refer to ingredients of these l.c.\ quantum groups $(M_1,\de_1)$, resp.\ $(M_2,\de_2)$. All notations
{\it without indices} refer to the {\it bicrossed product} quantum group $(M,\de)$. Finally, we will discuss an associated {\it double crossed
product} and then, the {\it index $\tw$} will be used.
\end{notation}

\section{Double crossed products: definition and multiplicative unitary} \label{sec.doubleone}

The matching $\tw$ of $(M_1,\deo)$ and $(M_2,\det)$ leads naturally to a coassociative comultiplication $\detw$ in the following definition.

\begin{definition}
Define the von Neumann algebra $M_\tw = M_1 \ot M_2$ and
$$\detw : \Mtw \recht \Mtw \ot \Mtw : \detw = (\io \ot \si\tw \ot
\io)(\deopo \ot \det) \; .$$
Then, $\detw$ is a co-associative comultiplication.
\end{definition}

One of our main objectives will be to prove that $(\Mtw,\detw)$ is always a l.c.\ quantum group. We will explicitly calculate all the ingredients of $(\Mtw,\detw)$ in terms of $(M_1,\de_1)$, $(M_2,\de_2)$ and the bicrossed product $(M,\de)$.

\begin{notation} \label{not.Z}
We write $Z=J \Jh (\Jh_1 J_1 \ot \Jh_2 J_2)$.
\end{notation}

\begin{lemma} \label{lemma.impltw}
We have $\tw(z) = Z z Z^*$ for all $z \in M_1 \ot M_2$.
\end{lemma}
\begin{proof}
For $x \in M_1$, we have
$$Z (x \ot 1) Z^* = \Jh U_\al (J_1 x J_1 \ot 1) U_\al^* \Jh = \Jh (J_1
x J_1 \ot 1) \Jh = U_\be (x \ot 1) U_\be^* = \be(x) = \tw(x \ot 1) \;
.$$
We analogously have $\tw(1 \ot x) = Z (1 \ot x)Z^*$ for $x \in M_2$.
\end{proof}

It is easy to define a unitary implementing the comultiplication
$\detw$, but it is less easy to prove that this unitary is
multiplicative.

\begin{definitiona} \label{def.multun}
Define on $H_1 \ot H_2 \ot H_1 \ot H_2$, the unitary $$\Wtw = (\Si
V_1^* \Si)_{13} \; Z^*_{34} W_{2,24} Z_{34} \; .$$ \eindeformule
\end{definitiona}

\begin{proposition}  \label{prop.multun}
The unitary $\Wtw$ is multiplicative and $\detw(z) = \Wtw^* (1 \ot 1
\ot z)\Wtw$ for all $z \in \Mtw$.

Moreover, we have $(\tw \ot \io \ot \io)(\Wtw) = Z^*_{34} W_{2,24}
Z_{34} \;(\Si V_1^* \Si)_{13} \; .$
\end{proposition}

\begin{proof}
Let $z \in \Mtw$. Then,
\begin{align*}
\Wtw^* (1 \ot 1 \ot z) \Wtw &= Z^*_{34} W^*_{2,24}
Z_{34} \; (\deopo \ot \io)(z)_{134} \; Z^*_{34} W_{2,24}
Z_{34} \\ &= Z^*_{34} \bigl( (\io \ot \io \ot \det)(\io \ot \tw)(\deopo
\ot \io)(z) \bigr)_{1324} Z_{34}= \detw(z) \; .
\end{align*}
We now prove that $\Wtw$ is a multiplicative unitary.

Because the bicrossed product is a l.c.\ quantum group, its
multiplicative unitary $W$ satisfies $\bigl( \Si (1 \ot J\Jh)W
\bigr)^3 \in \C$ (we can apply Proposition~6.9 of \cite{BS}).
In essentially the same way as in the first paragraph of the proof of Proposition~8.14 of
\cite{BS}, it follows that $\bigl( \Si (1 \ot \cU) \Wtw
\bigr)^3 \in \C$, where $\cU = (J_1 \Jh_1 \ot J_1 \Jh_1) Z$.
Define $N:= ({\Mh}'_1 \ot 1 \cup Z^*(1 \ot \Mh_2)Z )\dpr$. We claim
that $\cU N \cU^* \subset N'$. Clearly, $\cU Z^* (1 \ot \Mh_2)Z
\cU^*=1 \ot {\Mh}_2'$ and this commutes with ${\Mh}'_1 \ot 1$. To
prove that the commutator $[1 \ot {\Mh}_2',Z^* (1 \ot \Mh_2)Z]=\{0\}$,
it suffices to check that $[1 \ot \Mh_2, J \Jh (1 \ot \Mh_2) \Jh
J]=\{0\}$. We know this, because $1 \ot \Mh_2 \subset \Mh$, the dual
bicrossed product. Observing that, up to a scalar, $\cU = Z^* (J_1
\Jh_1 \ot J_1 \Jh_1)$, we check in an analogous way that $\cU
({\Mh}'_1 \ot 1) \cU^*  \subset N'$ and this proves our claim.
Because, $\bigl( \Si (1 \ot \cU) \Wtw \bigr)^3 \in \C$, we get
$$\Si \hat{V}_\tw^* \hat{W}_\tw^* (1 \ot \cU) \Wtw (\cU
\ot \cU) \in \C \; ,$$
where $\hat{W}_\tw = \Si \Wtw^* \Si$ and $\hat{V}_\tw = (\cU^* \ot 1)
\Wtw^* (\cU \ot 1)$. From our claim, it follows that $(1 \ot \cU) \Wtw
(\cU \ot \cU) \in \B(H_1 \ot H_2) \ot N'$ and hence, we get
$$\hat{W}_\tw^* (1 \ot 1 \ot z) \hat{W}_\tw = \hat{V}_\tw (z \ot 1 \ot
1) \hat{V}_\tw^* \quad\text{for all}\; z \in N \; .$$
So, for $x \in \Mht$,
\begin{align*}
\hat{W}_\tw^* (1 \ot 1 \ot Z^*(1 \ot x)Z) \hat{W}_\tw &=
\hat{V}_\tw (Z^*(1 \ot x)Z \ot 1 \ot 1) \hat{V}_\tw^* \\ &= (\cU^* \ot Z^*)
W^*_{2,24} (1 \ot \Jh_2 J_2 x J_2 \Jh_2 \ot 1 \ot 1) W_{2,24} (\cU \ot
Z) \\ &= (Z^* \ot Z^*)(\deht(x))_{24} (Z \ot Z) \; .
\end{align*}
We conclude that $\Wtw (Z^*(1 \ot x)Z \ot 1 \ot 1) \Wtw^* = (Z^* \ot
Z^*)(\dehopt(x))_{24} (Z \ot Z)$. It is easy to check that $\Wtw (y \ot 1
\ot 1 \ot 1) \Wtw^* = {{\deh}'_1}\Op(y)_{13}$ for all $y \in
{\Mh_1}'$, where ${\deh}'_1$ is the natural comultiplication on the commutant $\Mh_1'$ given by ${\deh}'_1(y)=
(\Jh_1 \ot \Jh_1)\deh_1(\Jh_1 y \Jh_1)(\Jh_1 \ot \Jh_1)$ for $y \in \Mh_1'$.
Then, we check that
$${\Wtw}_{3456} \;  {\Wtw}_{1234} \; {\Wtw}^*_{3456} = (\Si V_1^*
\Si)_{13}  (\Si V_1^* \Si)_{15} \; Z^*_{34} Z^*_{56} W_{2,24} W_{2,26}
Z_{34} Z_{56} = {\Wtw}_{1234} \; {\Wtw}_{1256} \; .$$ Hence, $\Wtw$ is
a multiplicative unitary. Then also $$(\detw \ot \io \ot \io)(\Wtw) =
{\Wtw}_{1256} \; {\Wtw}_{3456}$$ and a small calculation yields the
formula $(\tw \ot \io \ot \io)(\Wtw) = Z^*_{34} W_{2,24}
Z_{34} \;(\Si V_1^* \Si)_{13} \; .$
\end{proof}

In order to prove later that $(\Mtw,\detw)$ is a l.c.\ quantum group, we
need the following remarkable lemma.

\begin{lemma} \label{lemma.commutation}
The following holds.
\begin{enumerate}
\item $\tw$ is a $^*$-isomorphism,
\item $\tw(\tauo_{-t} \ot \taut_t) = (\tauo_{-t} \ot \taut_t) \tw$ for
  all $t\in \R$,
\item $\tw (R_1 \ot R_2) = (R_1 \ot R_2) \tw^{-1}$.
\end{enumerate}
If we write $\eta_t = \tauo_{-t} \ot \taut_t$, then the elements $(\io
\ot \io \ot \om)(\Wtw)$, $\om \in \B(H_1 \ot H_2)_*$ form a core for
$\eta_{-i/2}$ and $$(R_1 \ot R_2) \tw \eta_{-i/2} (\io
\ot \io \ot \om)(\Wtw) = (\io
\ot \io \ot \om)(\Wtw^*) \; .$$
\end{lemma}
\begin{proof}
If $\om = \om_{\xi,\eta}$, $\xi,\eta \in H_1 \ot H_2$, we take a basis
$(e_i)_{i \in I}$ for $H_1 \ot H_2$ and get
$$(\io \ot \io \ot \om)(\Wtw) = \sum_{i \in I} (\io \ot
\om_{e_i,\eta})((\Si V_1^* \Si)_{12}) \ot (\io \ot
\om_{\xi,e_i})(Z^*_{23} W_{2,13} Z_{23}) \; .$$
Using the closedness of $\eta_{-i/2}$, it is easy to conclude that $(\io \ot \io \ot \om)(\Wtw) \in
D(\eta_{-i/2})$ and
$$(R_1 \ot R_2) \eta_{-i/2} (\io \ot \io \ot \om)(\Wtw) = (\io \ot \io
\ot \om)\bigl( (\Si V_1 \Si)_{13} Z^*_{34} W_{2,24}^* Z_{34} \bigr) \;
,$$
for all $\om \in \B(H_1 \ot H_2)_*$. We claim that the elements $(\io
\ot \io \ot \om)(\Wtw)$ provide a core for $\eta_{-i/2}$. Denote by
$D_0$ the domain of the closure of the restriction of $\eta_{-i/2}$ to
these elements. We have to prove that $D_0 = D(\eta_{-i/2})$. Taking
$\om$ of the form $(x \ot 1) Z^* \om (1 \ot a)$ for $x \in \B(H_1)$
and $a \in M_2$, we observe that
$$z:= (\io \ot \io \ot \om)\bigl( (\Si V_1^* \Si)_{13} Z^*_{34} \bigl(
\al(a)(x \ot 1) \bigr)_{34} W_{2,24} \bigr) \in D_0$$
and
$$(R_1 \ot R_2) \eta_{-i/2} (z) = (\io \ot \io \ot \om)\bigl( (\Si V_1 \Si)_{13} Z^*_{34} \bigl(
\al(a)(x \ot 1) \bigr)_{34} W^*_{2,24} \bigr) \; .$$
Any element of $\B(H_1) \ot M_2$ can be approximated in the strong$^*$
topology by a bounded net of elements in $\al(M_2)(\B(H_1) \ot 1)$ and in
particular the element $1 \ot b$ for $b \in M_2$. So, we find that for
all $b \in M_2$ and $\om \in \B(H_1 \ot H_2)_*$,
$$z:= (\io \ot \io \ot \om)\bigl( (\Si V_1^* \Si)_{13} Z^*_{34} ((1
\ot b)W_2)_{24} \bigr) \in D_0$$
and
$$(R_1 \ot R_2) \eta_{-i/2} (z) = (\io \ot \io \ot \om)\bigl( (\Si V_1
\Si)_{13} Z^*_{34} ((1 \ot b) W^*_2)_{24} \bigr) \; .$$
Taking in the previous formula $\om$ of the form $(1 \ot x) \om$, with
$x \in {\Mh}'_2$ and using the fact that every element of $\B(H_2)$
can be approximated by an element in $M_2 {\Mh}'_2$, we arrive at
$$(\io \ot \om_1)(\Si V_1^* \Si) \ot (\io \ot \om_2)(W_2) \in D_0 \;
.$$
So, we have proven our claim. In particular, we get that the elements
$(\io \ot \io \ot \om)\bigl( Z^*_{34} W_{2,24}
Z_{34} \;(\Si V_1^* \Si)_{13} \bigr)$ are dense in $M_1 \ot
M_2$, because the adjoints of these elements are dense in the range of $(R_1 \ot R_2)\eta_{-i/2}$, which in its turn is dense in $M_1 \ot M_2$.
Because of Proposition~\ref{prop.multun}, they are in the image
of $\tw$ and hence, $\tw$ is a $^*$-isomorphism.

Combining with Proposition~\ref{prop.multun}, we get that
$$\tw^{-1} (R_1 \ot R_2) \eta_{-i/2} (\io \ot \io \ot \om)(\Wtw) =
(\io \ot \io \ot \om)(\Wtw^*)$$
and the elements $(\io \ot \io \ot \om)(\Wtw)$ form a core. Completely
analogously, one shows that
$$\tw^{-1} (R_1 \ot R_2) \eta_{i/2} (\io \ot \io \ot \om)(\Wtw^*) =
(\io \ot \io \ot \om)(\Wtw)$$ and now, the elements $(\io \ot \io \ot
\om)(\Wtw^*)$ form a core for $\eta_{i/2}$. It follows that, as unbounded mappings,
\begin{equation} \label{eq.twist}
\tw^{-1} (R_1 \ot R_2) \eta_{-i/2} = \bigl(\tw^{-1} (R_1 \ot R_2)
\eta_{i/2} \bigr)^{-1} = \eta_{-i/2} (R_1 \ot R_2) \tw \; .
\end{equation}
So, still as unbounded mappings, we have
$$\eta_{-i} = (\eta_{-i/2} (R_1 \ot R_2) \tw)(\tw^{-1} (R_1 \ot R_2)
\eta_{-i/2}) = (\tw^{-1} (R_1 \ot R_2)
\eta_{-i/2}) (\eta_{-i/2} (R_1 \ot R_2) \tw)$$
and we conclude that $(R_1 \ot R_2) \tw \eta_{-i} = \eta_{-i} (R_1 \ot
R_2) \tw$. The right hand side, is equal to $(R_1 \ot R_2) \eta_{-i}
\tw$ and we get $\tw \eta_{-i} = \eta_{-i} \tw$. From Sections 4.3 and
4.4 in \cite{haagerup},
it follows that $\tw (\tauo_{-t} \ot \taut_t) = (\tauo_{-t} \ot
\taut_t) \tw$ for all $t \in \R$. But then, also $\eta_{-i/2} \tw =
\tw \eta_{-i/2}$ and we conclude from Equation~\eqref{eq.twist} that
$\tw^{-1} (R_1 \ot R_2) = (R_1 \ot R_2) \tw$. Hence, we are done.
\end{proof}

\section{More about bicrossed products}

We emphasize again that all notations without indices refer to the bicrossed product quantum group $(M,\de)$ associated to our fixed matched pair $(M_1,\de_1)$ and $(M_2,\de_2)$, cf.\ Notation \ref{nota.crucial}.

We give an explicit description of all corepresentations of the bicrossed product $(M,\de)$ and deduce from this result some useful information on the modular elements
$\sde$ and $\sdeh$ of $(M,\de)$ and $(\Mh,\deh)$.

\begin{proposition}
Suppose that $X$ is a corepresentation of $(M,\de)$ on the Hilbert space $K$, i.e.\ $X \in M \ot \B(K)$ and $(\de \ot \io)(X) = X_{125} X_{345}$.
Then, there exists unique corepresentations $y$ of $(M_2,\de_2)$ and $z$ of $(\Mh_1,\deh_1)$ on $K$ such that
$$X=(\al \ot \io)(y) z_{13} \; .$$
Conversely, if $y$ and $z$ are corepresentations of $(M_2,\de_2)$ and
$(\Mh_1,\deh_1)$ respectively, on the Hilbert space $K$, the formula
$X=(\al \ot \io)(y) z_{13}$ defines a corepresentation $X$ of
$(M,\de)$ if and only if
$$z^*_{13} \; y_{23} \; Z_{12}^* z_{13} Z_{12} \in M_1' \ot \B(H_2 \ot
K) \; .$$
\end{proposition}

\begin{proof}
Because we have a morphism from $(M,\de)$ to $(\Mh_1,\deh_1)$, we can
define two actions of $(\Mh_1,\deh_1)$ on $M$. Explicitly, we have
that $\te : M \recht M \ot \Mh_1$ is a right action of
$(\Mh_1,\deh_1)$ while $\mu : M \recht \Mh_1 \ot M$ is a left action
and $\te$, $\mu$ are determined by
\begin{equation} \label{eq.muente}
(\te \ot \io \ot \io)(W) = W_{1245} (\io \ot \be)(\Wh_1)_{345}
\quad\text{and}\quad (\mu \ot \io \ot \io)(W) = (\io \ot
\be)(\Wh_1)_{145} W_{2345} \; .
\end{equation}
Because $\te$ is in fact the dual action on the crossed product $M =
M_1 \kruisje{\al} M_2$, we know that $M^\te = \al(M_2)$. Further, we
observe that $\mu = \si(R \ot \Rh_1)\te R$ and so, $M^\mu = \al(M_2)$
as well. From Equation
\eqref{eq.muente}, we conclude that $(\de \ot \io)\te = (\io \ot \io
\ot \te)\de$. So, we calculate
$$(\de \ot \io \ot \io) \bigl( X^*_{124} (\te \ot \io)(X) \bigr) = X^*_{346}
X^*_{126} (\io \ot \io \ot \te \ot \io)(X_{125} X_{345}) = 1 \ot 1 \ot  X^*_{124} (\te \ot \io)(X)$$ and conclude that there exists a unique $z \in
\Mh_1 \ot \B(K)$ such that $(\te \ot \io)(X) = X_{124} z_{34}$. Because applying $\te \ot \io \ot \io$ or $\io \ot \io \ot \deh_1 \ot \io$ to the
left hand side of this equation gives the same result, we get $(\deh_1 \ot \io)(z) = z_{13} z_{23}$.

Define $\gamma_1 : \Mh_1 \recht M : \gamma_1(a) = a \ot 1$. It is easy to check that $\te \gamma_1 = (\gamma_1 \ot \io)\deh_1$. If we apply $\io \ot
\io \ot \te \gamma_1 \ot \io$ to both sides of the equation $(\te \ot \io)(X) = X_{124} z_{34}$, we conclude that $(\te \ot \io)(z_{13}) =
z_{14}z_{34}$. But then, it follows that $(\te \ot \io)(X z^*_{13}) = (Xz^*_{13})_{124}$ and, because $M^\te=\al(M_2)$, we can take $y \in M_2 \ot
\B(K)$ such that $Xz^*_{13} = (\al \ot \io)(y)$. We conclude that $X = (\al \ot \io)(y) z_{13}$. We want to prove that $(\de_2 \ot \io)(y) = y_{13}
y_{23}$.

Because $X$ is a corepresentation and because of the formula for $(\te \ot \io)(X)$, we get $((\te \ot \io \ot \io)\de \ot \io)(X) = X_{126} z_{36}
X_{456}$. Using Equation \eqref{eq.muente}, we get that $(\te \ot \io \ot \io)\de = (\io \ot \io \ot \mu)\de$ and so, we conclude that $(\mu \ot
\io)(X) = z_{14} X_{234}$. Because $X = (\al \ot \io)(y) z_{13}$ and $\al(M_2) = M^\mu$, it follows that
\begin{equation} \label{eq.hulp}
(\mu \ot \io)(z_{13}) = (\al \ot \io)(y^*)_{234} z_{14} (\al \ot \io)(y)_{234} \; z_{24} \; .
\end{equation}
Next, we get from Equation \eqref{eq.muente} that $\mu(x) = (\io \ot
\be)(\Wh_1^*) (1 \ot x) (\io \ot \be)(\Wh_1)$ for all $x \in M$. In
particular, it follows that $\de \gamma_1 = (\gamma_1 \ot \io)\mu \gamma_1$.
Applying this to Equation \eqref{eq.hulp}, it follows that
$$(\de \ot \io)(z_{13}) = (\al \ot \io)(y^*)_{345} z_{15} (\al \ot
\io)(y)_{345} \; z_{35} \; .$$ Because $X$ is a corepresentation and $X = (\al \ot \io)(y) z_{13}$, we get
\begin{align*}
(\al \ot \io)(y)_{125} z_{15}  \; (\al \ot \io)(y)_{345} z_{35} &= X_{125} X_{345} = (\de \ot \io)(X) \\ &= ((\al \ot \al)\de_2 \ot \io)(y) \; (\al
\ot \io)(y^*)_{345} z_{15} (\al \ot \io)(y)_{345} \; z_{35} \; .
\end{align*}
It follows that $((\al \ot \al)\de_2 \ot \io)(y) = (\al \ot \io)(y)_{125} (\al \ot \io)(y)_{345}$, which finally yields $(\de_2 \ot \io)(y) = y_{13}
y_{23}$.

The uniqueness statement is obvious from the given construction. To
prove the converse statement, the formula $z^*_{13} \; y_{23} \;
Z_{12}^* z_{13} Z_{12} \in M_1' \ot \B(H_2 \ot K)$ is equivalent with
$$z^*_{24} \; y_{34} \; Z_{23}^* z_{24} Z_{23} \quad\text{commutes
  with}\quad \Wh_{1,12} \; ,$$
which is equivalent with
$$z_{24} \Wh_{1,12}^* z_{24}^* \;  y_{34} \; Z_{23}^* z_{24} Z_{23} \;
\Wh_{1,12} = y_{34} \; Z_{23}^* z_{24} Z_{23} \; .$$
Using that $z$ is a corepresentation of $(\Mh_1,\deh_1)$, the last
formula is equivalent with
$$\Wh_{1,12}^* \; Z_{23}^* z_{24} Z_{23} \; \Wh_{1,12} = y^*_{34}
z_{14} y_{34} \; Z_{23}^* z_{24} Z_{23} \; .$$
Applying $\Ad Z_{23}$ and using the fact that $\de(a \ot 1) = (\io \ot
\be)(\Wh_1^*)_{134} (1 \ot 1 \ot a \ot 1) (\io \ot \be)(\Wh_1)_{134}$
for $a \in \Mh_1$, the last formula is equivalent with
$$(\de \ot \io)(z_{13}) = (\al \ot \io)(y^*)_{345} z_{15} (\al \ot
\io)(y)_{345} \; z_{35} \; .$$
Because $y$ is a corepresentation of $(M_2,\de_2)$ and $\de \al = (\al
\ot \al)\de_2$, it is clear that the final formula is equivalent with
$X = (\al \ot \io)(y) z_{13}$ being a corepresentation of $(M,\de)$.
\end{proof}

\begin{proposition} \label{prop.rho}
There exist strictly positive, self-adjoint operators $\rho_1,\rho_2$ affiliated
with $M_1,M_2$ such that
\begin{itemize}
\item $\sde^{it} = \al(\rho_2^{it})(\sdeh_1^{it} \ot 1)$ and both
  factors commute;
\item $\sdeh^{it} = \be(\rho_1^{it})(1 \ot \sdeh_2^{it})$ and both
  factors commute;
\item $\deo(\rho_1) = \rho_1 \ot \rho_1$ and $\det(\rho_2) = \rho_2
  \ot \rho_2$;
\item $\nu = \frac{\nu_2}{\nu_1}$;
\item we have \vspace{-5.1ex}
\begin{align*}
J \sde^{it} J &= \Jh_1 \sdeh_1^{it} \Jh_1  \ot \rho_2^{-it}
\sde_2^{it} J_2 \sde_2^{it} J_2 \; , \\
\Jh \sdeh^{it} \Jh &= \rho_1^{-it} \sde_1^{it} J_1 \sde_1^{it} J_1 \ot
\Jh_2 \sdeh_2^{it} \Jh_2 \; , \\
[D \vfi_2 \na \al : D \vfi_2]_t &= \tw(\rho_1^{-it/2} \ot
\rho_2^{-it/2}) (\rho_1^{-it/2} \sde_1^{it}  \ot \rho_2^{it/2} ) \; ,
\end{align*}
where we used the Radon-Nikodym derivative introduced in Definition~\ref{def.radnik};
\item $\tw$ leaves invariant $\rho_1^{-it/2} \sde_1^{it} \ot
  \rho_2^{it/2} \sde_2^{-it}$.
\end{itemize}
There exist positive numbers $\sla_{1,2} > 0$ such that
$\si^1_t(\rho_1^{is}) = \sla_1^{ist} \rho_1^{is}$ and analogously for
$\sla_2$ and they satisfy $\frac{\sla_2}{\sla_1}=\frac{\nu_2}{\nu_1}$.
\end{proposition}
\begin{proof}
Because $\sde^{it}$ is a one-dimensional corepresentation of
$(M,\de)$, there exist unique one-dimensional corepresentations (i.e.\
group-like elements) $v_t \in
\Mh_1$ and $u_t \in M_2$ such that $\sde^{it} = (v_t \ot 1) \al(u_t)$.
Because $\al$ and $\be$ are morphisms we know that
\begin{align*}
& \Ad(\sde^{it}) \na \al = \al \na \Ad(\sde_2^{it}) \; , \quad
\Ad(\sdeh^{it}) \na \al = \al \na \Ad(\sdeh_2^{it}) \quad\text{on}\;
M_2 \; , \\ & \Ad(\sde^{it}) \na \be = \be \na \Ad(\sdeh_1^{it}) \; , \quad
\Ad(\sdeh^{it}) \na \be = \be \na \Ad(\sde_1^{it}) \quad\text{on}\;
M_1 \; .
\end{align*}
Then,
\begin{align*}
\deop(v_t \ot 1) &= (\be \ot \io)(W_1) (v_t \ot 1 \ot 1) (\be \ot
\io)(W_1^*) \ot 1 \\ &= \bigl( (\be \ot \io)(W_1) (\sde^{it} \ot 1) (\be
\ot \io)(W_1^*) \ot 1 \bigr) \; (\al(u_t^*) \ot 1 \ot 1) \\ &= \bigl( (\be \ot \io)(W_1)
(\be \na \Ad (\sdeh_1^{it}) \ot \io)(W_1^*) (v_t \ot 1 \ot 1) \bigr)
\ot 1 = v_t \ot 1 \ot \sdeh_1^{it} \ot 1 \; .
\end{align*}
Hence, $\deop(\sde^{it}) = \deop( (v_t \ot 1) \al(u_t)) = (v_t \ot 1)
\al(u_t) \ot (\sdeh_1^{it} \ot 1) \al(u_t)$. Because $\sde^{it}$ is
group-like, we conclude that $v_t = \sdeh_1^{it}$.

At the end of this proof, we will argue that in a l.c.\ quantum group,
the modular element commutes with any group-like unitary. So,
$\sde^{it}$ commutes with $\al(u_s)$ for all $s,t \in \R$. As
$\sde^{it}$ commutes with $\sde^{is}$, we get that $\sde^{it}$
commutes with $\sdeh_1^{is}$ for all $s,t \in \R$. This means that
$((\sdeh_1^{-it} \ot 1) \sde^{it})_t$ is a one-parameter group of
unitaries. Then, there exists a unique strictly positive, self-adjoint
operator $\rho_2$ affiliated with $M_2$ such that $\sde^{it} =
(\sdeh_1^{it} \ot 1) \al(\rho_2^{it})$. By symmetry, we find $\rho_1$
and we have proven the first three statements of the proposition.

Let $a \in \cN_{\vfih_1}$ and $x \in \cN_{\vfi_2}$. Then, we have
$$(a \ot 1) \al(x) \sde^{it} = (a \ot 1) \sde^{it} \al(\sde_2^{-it} x
\sde_2^{it}) = (a \sdeh_1^{it} \ot 1) \al(\rho_2^{it} \sde_2^{-it} x
\sde_2^{it}) \; .$$
We conclude that
$$\la\bigl((a \ot 1) \al(x) \sde^{it}\bigr) =
\bigl(\frac{\nu_2}{\nu_1})^{-t/2} (\Jh_1 \sdeh_1^{-it} \Jh_1 \ot
\rho_2^{it} \sde_2^{-it} J_2 \sde_2^{-it} J_2) \la\bigl( (a \ot 1)
\al(x) \bigr) \; .$$
On the other hand, we know that
$$\la\bigl((a \ot 1) \al(x) \sde^{it}\bigr) = \nu^{-t/2} J \sde^{-it}
J \la\bigl( (a \ot 1) \al(x) \bigr) \; .$$
So, $\nu = \frac{\nu_2}{\nu_1}$ and $J \sde^{it} J =
\Jh_1 \sdeh_1^{it} \Jh_1 \ot \rho_2^{-it} \sde_2^{it} J_2 \sde_2^{it}
J_2$. We find the formula for $\Jh \sdeh^{it} \Jh$ in an analogous
way.

From \cite{VVD}, Proposition 2.4, we know that
$$\nab^{it} = \sdeh^{-it/2} \sde^{-it/2} J \sde^{-it/2} J \Jh
\sdeh^{it/2} \Jh \; .$$
We can take positive numbers $\sla_{1,2} > 0$ such that
$\si^1_t(\rho_1^{is}) = \sla_1^{ist} \rho_1^{is}$ and analogously for
$\sla_2$ (see the last paragraph of this proof). Because in any l.c.\ quantum group, $\Ad \sdeh^{it} =
\tau_{-t} \Ad(\sde^{-it}) \si_{-t}$ on $M$, we easily find the
required formula for $[D \vfi_2 \na \al : D \vfi_2]_t := \nab^{it}
(\nabh_1^{-it} \ot \nab_2^{-it})$. We also get that $\sdeh^{-is}
\al(\rho_2^{it}) \sdeh^{is} = \sla_2^{ist} \al(\rho_2^{it})$ and
$\sde^{it} \be(\rho_1^{is}) \sde^{-it} = \sla_1^{-ist}
\be(\rho_1^{is})$. A short calculation yields that $\sde^{it}
\sdeh^{is} = \bigl(\frac{\sla_2}{\sla_1})^{ist} \sdeh^{is}
\sde^{it}$. But, from the general theory, we know that
$\sde^{it}
\sdeh^{is} = \nu^{ist} \sdeh^{is}
\sde^{it}$, giving us the equality $\frac{\sla_2}{\sla_1}=\frac{\nu_2}{\nu_1}$.

Consider $\detop$ as an action of $(M_2,\detop)$ on $M_2$. Then, $[D
\vfi_2 \na \detop : D \vfi_2]_t = \sde_2^{-it} \ot 1$. Because $(\io
\ot \detop)\al = (\tw \si \ot \io)(\io \ot \al)\detop$ and because, by
Lemma~\ref{lemma.commutation}, $\tw(\tauo_{-t} \ot \taut_t) =
(\tauo_{-t} \ot \taut_t) \tw$, we can apply
Lemma~\ref{lemma.formula} of our Appendix to conclude that
$$(\io \ot \detop)([D \vfi_2 \na \al : D \vfi_2 ]_t) = (\tw \si \ot
\io)(\sde_2^{-it} \ot [D \vfi_2 \na \al : D \vfi_2 ]_t) (1 \ot
\sde_2^{it} \ot 1) \; .$$
Filling in the already obtained formula for $[D \vfi_2 \na \al : D
\vfi_2 ]_t$ and carefully rewriting, we arrive at the final statement
of the proposition.

As we promised above, we consider now a group-like unitary $u \in M$
for some l.c.\ quantum group $(M,\de)$. By the uniqueness of the Haar
measure, we find $\sla > 0$ such that $\vfi(uxu^*) = \sla \vfi(x)$ for
all $x \in M^+$. It follows that $\si_t(u) = \sla^{it} u$ for all $t
\in \R$. Because a group-like unitary is a one-dimensional corepresentation,
we have $u \in D(S)$ and $S(u) = u^*$. But, also $u^*$ is a group-like
element and so $S^2(u) = u$. It follows that $\tau_t(u) = u$ for all
$t \in \R$ and $R(u) = u^*$. Then $\si'_t(u) = R \si_{-t} R(u) =
\sla^{it} u$, which gives $\si'_t \si_{-t} (u) = u$ and so, $u$
commutes with $\sde^{it}$ for all $t \in \R$.
\end{proof}

\section{Double crossed products as locally compact quantum groups} \label{sec.doublelc}

In this section, we compute the invariant weights on $(\Mtw,\detw)$
and obtain that $(\Mtw,\detw)$ is a l.c. quantum group. We compute all
the operators associated with $(\Mtw,\detw)$.

\begin{notation}
We use $\rho_1,\rho_2$ to denote the unbounded operators introduced in
Proposition \ref{prop.rho}. We define strictly positive, self-adjoint
operators $k_1,k_2$ such that
$$k_1^{it} = \rho_1^{it} \sde_1^{-it} \quad\text{and}\quad k_2^{it} =
\rho_2^{it} \sde_2^{-it} \; .$$
\end{notation}

In order to obtain invariant weights for $(\Mtw,\detw)$, we need the following relative invariance result of the weight $\psi_2$ under the action $\al$.

\begin{proposition} \label{prop.afgeleide}
For all $x \in \cN_{\psi_2}$, $\xi \in H_1$ and $\eta \in
D(k_1^{1/2})$, we have $(\om_{\xi,\eta} \ot \io)(U^*_\be \al(x) ) \in
\cN_{\psi_2}$ and
$$\Gamma_2\bigl(  (\om_{\xi,\eta} \ot \io)(U^*_\be \al(x) ) \bigr) =
(\om_{\xi,k_1^{1/2} \eta} \ot \io)(U_\be^* U_\al) \Gamma_2(x) \; .$$
On the other hand, for $x \in \cN_{\psi_1}$, $\xi \in H_2$ and $\eta
\in D(k_2^{1/2})$, we have
$$\Gamma_1\bigl( (\io \ot \om_{\xi,\eta})(U_\al^* \be(x)) \bigr) =
(\io \ot \om_{\xi,k_2^{1/2} \eta})(U_\al^* U_\be) \Ga_1(x) \; .$$
\end{proposition}
\begin{proof}
Consider the dual weight $\psitil_2$ on $M$ with a canonical
GNS-map $\Gatil$ determined by
$$\Gatil \bigl( (a \ot 1) \al(x) \bigr) = \lah_1(a) \ot \Gamma_2(x)
\quad\text{for all}\quad a \in \cN_{\vfih_1}, x \in \cN_{\psi_2} \;
.$$
We denote by $\nabtil$ the modular operator of $\psitil_2$ in this
GNS-construction. Observe that $\psi_2 = (\vfi_2)_{\sde_2}$ and hence,
$\psitil_2 = \vfi_{\al(\sde_2)}$.
It then follows
from \cite{V2}, Proposition 2.5 that
\begin{align*}
\nabtil^{it} &= J \al(\sde_2^{it}) J \; \al(\sde_2^{it}) \; \nab^{it}
\; , \\
\tilde{S} &= \nu_2^{i/4} J \nabtil^{1/2} \; ,
\end{align*}
where $\tilde{S}$ is determined by $\tilde{S} \Gatil(z) =
\Gatil(z^*)$. From Proposition \ref{prop.rho}, it follows that
\begin{align*}
\nabtil^{it} &= \al(\sde_2^{it}) (1 \ot J_2 \sde_2^{it} J_2)
\nab^{it} = \al(\sde_2^{it}) (1 \ot J_2 \sde_2^{it} J_2) [D \vfi_2 \na
\al : D \vfi_2]_t (\nabh_1^{it} \ot \nab_2^{it}) \\ &= \tw(\rho_1^{-it/2}
\ot \rho_2^{-it/2} \sde_2^{it}) (\rho_1^{-it/2} \sde_1^{it}
\ot \rho_2^{it/2} \sde_2^{-it}) (\nabh_1^{it} \ot \nabp_2^{it}) =
\be(k_1^{-it}) (\nabh_1^{it} \ot \nabp_2^{it}) \; ,
\end{align*}
where $\nabp_2$ denotes the modular operator of $\psi_2$.

We know that the elements $\al(x^*) (\mu \ot \Ga_2(y))$, $x,y \in
\cN_{\psi_2}$, $\mu \in D(\nabh_1^{1/2})$ span a core for $\Stil$ and
$$\Stil \al(x^*) (\mu \ot \Ga_2(y)) = \al(y^*) (\Jh_1 \nabh_1^{1/2}
\mu \ot \Ga_2(x)) \; .$$
It follows that
$$\nabtil^{1/2} \al(x^*) (\mu \ot \Ga_2(y)) = \nu_2^{i/4} (1 \ot J_2
y^* J_2) U_\al (\nabh_1^{1/2} \mu \ot J_2 \Ga_2(x)) \; .$$
Because the product of $\be(k_1^{-it})$ and $\nabh_1^{it} \ot
\nabp_2^{it}$ yields the one-parameter group $\nabtil^{it}$, we get
that $\be(k_1)$ and $\nabh_1 \ot \nabp_2$ commute strongly.
If we assume now that $x$ is analytic w.r.t.\ the modular group $(\si^{2 \, \prime}_t)$ of $\psi_2$ and $y \in \cN_{\psi_2} \cap \cN_{\psi_2}^*$,
we
conclude that the left hand side of the previous equality is equal to
$$\nu_2^{i/4} \be(k_1^{1/2}) \al\bigl( \si^{2 \, \prime}_{-i/2}(x^*)
\bigl) (\nabh_1^{1/2} \mu \ot J_2 \Ga_2(y^*)) \; .$$
By a typical density argument, we may conclude that for $x \in
\cN_{\psi_2}$, analytic w.r.t.\ $(\si^{2 \, \prime}_t)$, $y \in
\cN_{\psi_2}$ and $\xi \in H_1$,
$$\be(k_1^{-1/2}) \al\bigl(\si^{2 \, \prime}_{-i/2}(x^*)\bigr) (\xi
\ot J_2 \Ga_2(y)) = (1 \ot J_2 y J_2) U_\al (\xi \ot J_2 \Ga_2(x)) \;
.$$
Applying $U_\be^* \in B(H_1) \ot M_2$, it follows that for $x,y,\xi$
as above and $\eta \in D(k_1^{1/2})$,
$$(\om_{\xi,\eta} \ot \io)\bigl(U_\be^* \al(\si^{2 \,
  \prime}_{-i/2}(x^*)) \bigr) J_2 \Ga_2(y) = J_2 y J_2
(\om_{\xi,k_1^{1/2} \eta} \ot \io)(U_\be^* U_\al) J_2 \Ga_2(x) \; .$$
Using $J_2 \Ga_2(z) = \Ga_2(\si^{2 \,
  \prime}_{-i/2}(z^*))$, we get, for $x,y \in \cN_{\psi_2}$, $y$ analytic w.r.t.\ $(\si^{2 \,
  \prime}_t)$, $\xi \in H_1$ and $\eta \in D(k_1^{1/2})$,
$$\Ga_2 \bigl( (\om_{\xi,\eta} \ot \io)(U_\be^* \al(x)) y \bigr) = J_2
\si^{2 \, \prime}_{i/2}(y)^* J_2 (\om_{\xi,k^{1/2} \eta} \ot
\io)(U_\be^* U_\al) \Ga_2(x) \; .$$
From this formula, our first result follows. The second result is
analogous.
\end{proof}

\begin{theorem} \label{thm.double}
The pair $(\Mtw,\detw)$ is a l.c.\ quantum group. The weight $\vfitw:=\psi_1
\ot (\vfi_2)_{k_2}$ is left invariant and
$$\Wtw= (\Si
V_1^* \Si)_{13} \; Z^*_{34} W_{2,24} Z_{34} =(\Si V_1^* \Si)_{13} U_{\be, 32}^* \; (\al \ot
\io)(W_2)_{324} $$
is the left regular
representation when we use the canonical GNS-construction for
$\vfitw$.

The canonical right invariant weight $\psitw$ is given by
$(\vfi_1)_{k_1} \ot \psi_2$  and the
right regular representation by
$$\Vtw = Z^*_{12} \Wh_{1,13} Z_{12} V_{2,24}
= (\io \ot \be)(\Wh_1)_{132} U_{\al,32} \; V_{2,24} \; .$$
Further, we have
\begin{alignat*}{2}
\Rtw &= \tw^{-1} (R_1 \ot R_2) = (R_1 \ot R_2) \tw \; , & \qquad \qquad
\nutw &= \frac{\nu_2}{\nu_1} \; , \\
\tautw_t &= \tauo_{-t} \ot \taut_t \; , & \qquad\qquad
P_\tw &= P_1^{-1} \ot P_2 \; , \\
\Jtw &= \sla_1^{i/4} (J_1 \ot J_2) \; , & \qquad \qquad
\sdehtw^{it} &= (\Jh_1 \sdeh_1^{-it} \Jh_1 \ot 1) Z^*(1 \ot \sdeh_2^{it})Z  \; , \\
\Jhtw &= \sla_1^{i/4} (\Jh_1 \ot \Jh_2) J \Jh (\Jh_1 J_1 \ot \Jh_2 J_2)
\; , & \qquad \qquad
\Jhtw \sdehtw^{it} \Jhtw &= (1 \ot \Jh_2 \sdeh_2^{it} \Jh_2) Z^*(\sdeh_1^{-it} \ot 1)Z \; , \\
\sdetw^{it} &= \rho_1^{it} \sde_1^{-2it} \ot \rho_2^{-it}
\sde_2^{2it} \; , & &
\end{alignat*}
where $\sla_1 > 0$ is defined in Proposition~\ref{prop.rho}.
The dual of $(\Mtw,\detw)$ is given by
$$\Mhtw = \bigl(\Mh_1' \ot 1 \cup Z^* (1 \ot \Mh_2)Z \bigr)\dpr \qquad\text{and}\qquad \Mh_\tw' = \bigl( Z^*(\Mh_1 \ot 1)Z \cup 1 \ot \Mh_2'\bigr)\dpr \; ,$$
and in this picture, $(\Mh_1',\deh_1')$ and $(\Mh_2,\deh_2)$ are closed quantum subgroups of $(\Mhtw,\dehtw)$.

Finally, we have the following commutation relations
\begin{align}
J \Jh (J_1 \ot J_2) &= (\frac{\nu_1}{\nu_2})^{i/4} (J_1 \ot J_2) J
\Jh \; , \notag \\
(\io \ot \tw)(U_{\be,13} V_{1,12}) &= V_{1,12} U_{\be,13} \; , \label{eq.commone}\\
(\tw \ot \io)(U_{\al,13} (\Si V_2 \Si)_{23} ) &= (\Si V_2 \Si)_{23}
U_{\al,13} \; . \label{eq.commtwo}
\end{align}
\end{theorem}

Recall that the notations $Z$ and $\Wtw$ are introduced in
Notation~\ref{not.Z} and Proposition~\ref{prop.multun}.

\begin{proof}
Define $\vfitw:=\psi_1 \ot (\vfi_2)_{k_2}$ and denote by $\latw$ its
canonical GNS-map. Let $x \in \cN_{\psi_1}$,
$y,z \in \cN_{(\vfi_2)_{k_2}}$ where $z$ is analytic w.r.t.\ the
modular group of $(\vfi_2)_{k_2}$. Take $\xi,\eta \in H_1$ and
$\mu,\rho \in H_2$. Then,
$$(\om_{\xi \ot \mu, \eta \ot \rho} \ot \io \ot \io)\detw(x \ot y) (1
\ot z) = (\io \ot \om_{\mu,\rho} \ot \io) \bigl( (U_\al \ot 1) (1 \ot \det(y))
\bigl( U_\al^* \be((\io \ot \om_{\xi,\eta})\deo(x))  \ot z \bigr)
\bigr) \; .$$
Let $\chi_{[1/n,n]}$ be the characteristic function of the interval
$[1/n,n]$ and define $P_n = \chi_{[1/n,n]}(k_2)$. Let $(e_i)$ be an
orthonormal basis for $H_2$. Consider the following element in $M_1
\ot M_2$:
\begin{align*}
r & := (\io \ot \om_{\mu,\rho} \ot \io) \bigl( (U_\al \ot 1) (1 \ot
\det(y)) (1 \ot P_n \ot 1)
\bigl( U_\al^* \be((\io \ot \om_{\xi,\eta})\deo(x))  \ot z \bigr)
\bigr) \\ &= \sum_i (\io \ot \om_{e_i,\rho} \ot \io)\bigl( (U_\al \ot 1) (1 \ot
\det(y)) \bigr) \; \bigl( (\io \ot \om_{\mu,P_n e_i})
\bigl( U_\al^* \be((\io \ot \om_{\xi,\eta})\deo(x))\bigr) \ot z \bigr)
\; .
\end{align*}
Denote every individual term of this sum by $r_i$. Because of
Proposition~\ref{prop.afgeleide}, every $r_i$ belongs to
$\cN_{\vfitw}$ and
$$\latw(r_i) = (\io \ot \om_{e_i,\rho} \ot \io)\bigl( (U_\al \ot 1) (1 \ot
\det(y)) \bigr) \; \bigl( (\io \ot \om_{\mu,k_2^{1/2} P_n
  e_i})(U^*_\al U_\be) (\io \ot \om_{\xi,\eta})(V_1) \Ga_1(x) \ot
\la_{(\vfi_2)_{k_2}}(z) \bigr) \; .$$
This expression is summable over $i$ and we conclude that $r \in
\cN_{\vfitw}$ and
$$\latw(r) = (\io \ot \om_{\mu,\rho} \ot \io)\bigl( (\al \ot
\io)(\de_2(y)((k_2^{1/2} P_n) \ot 1) ) (U_\be \ot 1) \bigr) \; \bigl(
(\io \ot \om_{\xi,\eta})(V_1) \Ga_1(x) \ot
\la_{(\vfi_2)_{k_2}}(z) \bigr) \; .$$
Write $v:= (\io \ot \om_{\xi,\eta})(V_1) \Ga_1(x)$ and take $w \in H_1$.
Then, denoting by $(\sitil_t)$ the modular group of $(\vfi_2)_{k_2}$,
\begin{align*}
(\te_w^* \ot \io)\latw(r) & =
\la_{(\vfi_2)_{k_2}} \bigl( (\om_{U_\be (v \ot \mu), w \ot \rho} \na
\al \ot \io)(\det(y) ((k_2^{1/2} P_n) \ot 1)) \; z \bigr) \\ &= J_2
\sitil_{i/2}(z)^* J_2 \; (\om_{\al(P_n) U_\be (v \ot \mu), w \ot \rho} \na
\al \ot \io) (W_2^*) \; \la_{(\vfi_2)_{k_2}}(y) \; ,
\end{align*}
where we used that $k_2$ is a group-like.
Hence,
$$\latw(r) = (1 \ot J_2
\sitil_{i/2}(z)^* J_2) \; (\io \ot \om_{\mu,\rho} \ot \io) \bigl(
(\al \ot \io)(W_2^*) (\al(P_n)U_\be \ot 1) \bigr) \; \bigl( (\io \ot
\om_{\xi,\eta})(V_1) \Ga_1(x) \ot \la_{(\vfi_2)_{k_2}}(y) \bigr) \;
.$$
Taking limits over $n \recht \infty$ and $z \recht 1$, we get that
$$(\om_{\xi \ot \mu, \eta \ot \rho} \ot \io \ot \io)\detw(x \ot y) \in
\cN_{\vfitw}$$
and
$$\latw \bigl( (\om_{\xi \ot \mu, \eta \ot \rho} \ot \io \ot
\io)\detw(x \ot y) \bigr) =
(\om_{\xi \ot \mu, \eta \ot \rho} \ot \io \ot \io) \bigl( (\al \ot
\io)(W_2^*)_{324} \; U_{\be, 32} (\Si V_1 \Si)_{13} \bigr) \latw(x \ot
y) \; .$$
It follows that $\vfitw$ is left invariant and that its associated left
regular representation is given by
$$\cWtw := (\Si V_1^* \Si)_{13} U_{\be, 32}^* \; (\al \ot
\io)(W_2)_{324} \; .$$
Define $\Rtw:= (R_1 \ot R_2) \tw$. From Lemma~\ref{lemma.commutation},
we know that $\Rtw = \tw^{-1} (R_1 \ot R_2)$ as well and then, it is
easy to verify that $\detw \Rtw = (\Rtw \ot \Rtw)\detwop$.

So, $\vfitw \Rtw$ is a right invariant weight and $(\Mtw,\Rtw)$ is a
l.c.\ quantum group.
It follows that $\cWtw$ is a multiplicative unitary and we want to
prove that $\cWtw=\Wtw$.

An analogous computation as above shows that $(\vfi_1)_{k_1} \ot
\psi_2$ is a right invariant weight with associated right regular
representation
$$\Vtw := (\io \ot \be)(\Wh_1)_{132} U_{\al,32} \; V_{2,24} \; .$$
Observe that
\begin{align*}
(\detw \ot \io \ot \io)(\cWtw) &= \si_{23} \tw_{23} \bigl(
(\Si V_1^* \Si)_{15} (\Si V_1^* \Si)_{25} U^*_{\be, 53} U^*_{\be, 54}
((\io \ot \det)\al \ot \io)(W_2)_{5346} \bigr) \\
&= \si_{23} \tw_{23} \bigl((\Si V_1^* \Si)_{15} (\Si V_1^* \Si)_{25}
U^*_{\be, 53} U^*_{\be, 54} (\tw_{13}(\al \ot \io \ot \io)(W_{2,13}
W_{2,23}))_{5346} \bigr) \\
&= (\Si V_1^* \Si)_{15} ((\io \ot \tw)(V^*_{1,12} U^*_{\be,13}))_{532}
((\io \ot \al)\al \ot \io)(W_2)_{5326} U^*_{\be,54} (\al \ot
\io)(W_2)_{546} \; .
\end{align*}
On the other hand, we know that
$$(\detw \ot \io \ot \io)(\cWtw) = (\cWtw)_{1256} (\cWtw)_{3456} =
(\Si V_1^* \Si)_{15} U_{\be, 52}^* \; (\al \ot
\io)(W_2)_{526} \; (\Si V_1^* \Si)_{35} U_{\be, 54}^* \; (\al \ot
\io)(W_2)_{546} \; .$$
We conclude that
$$(\io \ot \tw)(U_{\be,13} V_{1,12}) = V_{1,12} U_{\be,13} \; .$$
Making an analogous computation with $\Vtw$, we arrive at the
commutation relations \eqref{eq.commone} and \eqref{eq.commtwo}.

We now claim that $$Z^*_{23} W^*_{2,13} Z_{23} = (\al \ot
\io)(W_2^*)_{213} U_{\be, 21} \; ,$$ yielding the equality $\cWtw =
\Wtw$. Analogously, one proves the second formula for the right
regular representation: $\Vtw = Z^*_{12} \Wh_{1,13} Z_{12} V_{2,24}$. To prove the
claim, we make the following calculation, using the notation $u_i =
J_i \Jh_i$:
\begin{align*}
Z^*_{23} W^*_{2,13} Z_{23} &= (\Jh_2 \ot J_1 \ot \Jh_2) U^*_{\al,23}
U_{\be,23} (\Si V_2 \Si)_{13} U^*_{\be,23} U_{\al,23} (\Jh_2 \ot J_1 \ot \Jh_2)
 \\ &=  (\Jh_2 \ot J_1 \ot \Jh_2) U^*_{\al,23} U_{\be,23}
(\io \ot \deopt)(U_\be^*)_{213} (\Si V_2 \Si)_{13}
U_{\al,23} (\Jh_2 \ot J_1 \ot \Jh_2)
 \\ &=  (\Jh_2 \ot J_1 \ot \Jh_2) U^*_{\al,23} U_{\be,21}^*
(\Si V_2 \Si)_{13} U_{\al,23} (\Jh_2 \ot J_1 \ot \Jh_2) \; .
\end{align*}
Observe that
$$Z^* = \nu_1^{-i/4} \Util_\al^* U_\be^* \quad\text{where}\quad \Util_\al = (J_1
\ot \Jh_2) U_\al (J_1 \ot \Jh_2) \in M_1' \ot B(H_2) \; .$$
We combine Lemma \ref{lemma.impltw} and Equation \eqref{eq.commtwo}
with the calculation above, to get
\begin{align*}
Z^*_{23} W^*_{2,13} Z_{23} &= \nu_1^{-i/4} (\Jh_2 \ot J_1 \ot \Jh_2) U^*_{\al,23}
\Util_{\al,21} Z^*_{21} (\Si V_2 \Si)_{13} U_{\al,23} (\Jh_2 \ot J_1
\ot \Jh_2) \\ &= \nu_1^{-i/4} (\Jh_2 \ot J_1 \ot \Jh_2) U^*_{\al,23}
\Util_{\al,21} U_{\al,23} (\Si V_2 \Si)_{13} Z^*_{21} (\Jh_2 \ot J_1 \ot \Jh_2)
\\ &=  \nu_1^{-i/4} (\Jh_2 \ot J_1 \ot \Jh_2) \Util_{\al,21}
(\Si V_2 \Si)_{13} Z^*_{21} (\Jh_2 \ot J_1 \ot \Jh_2) \\
&= U_{\al,21} W^*_{2,13} U^*_{\al,21} U_{\be,21} = (\al \ot
\io)(W^*_2)_{213} U_{\be,21} \; .
\end{align*}
This proves our claim.

Defining $\tautw_t:=\tauo_{-t} \ot \taut_t$, it follows from the
proof of Lemma \ref{lemma.commutation} that for all $\om \in B(H_1 \ot
H_2)_*$, $$\Rtw \tautw_{-i/2} (\io \ot \io \ot \om)(\Wtw) =
(\io \ot \io \ot \om)(\Wtw^*) \; .$$ From Lemma
\ref{lemma.commutation}, we also know that $\Rtw$ and $\tautw_t$
commute. So, we can conclude that $\Rtw$ is the unitary antipode of
$(\Mtw,\detw)$ and $(\tautw_t)$ is its scaling group.

Taking the $it$-th power of the Radon-Nikodym derivative of
$(\vfi_1)_{k_1} \ot \psi_2$ w.r.t.\ $\vfitw$, we get
$\rho_1^{it} \sde_1^{-2it} \ot \rho_2^{-it} \sde_2^{2it}$. These
unitaries are group-like because of Lemma \ref{lemma.commutation}. By
uniqueness of the Haar weights, $\vfitw \Rtw$ is proportional to
$(\vfi_1)_{k_1} \ot \psi_2$ and also the Radon-Nikodym derivative of
$\vfitw \Rtw$ w.r.t.\ $\vfitw$ is group-like. Because a non-trivial multiple of a
group-like is no longer group-like, we conclude that $\vfitw \Rtw =
(\vfi_1)_{k_1} \ot \psi_2$. So, $(\vfi_1)_{k_1} \ot \psi_2$ is indeed
the canonical right invariant weight of $(\Mtw,\detw)$.

We also conclude that $\sdetw^{\it}=\rho_1^{it} \sde_1^{-2it} \ot \rho_2^{-it}
\sde_2^{2it}$. Because clearly $\vfitw \tautw_t = (\nu_1 / \nu_2)^t
\vfitw$, we get that $\nutw = \frac{\nu_2}{\nu_1}$.

From \cite{V2}, Proposition 2.5 and from Proposition \ref{prop.rho}
above, we know that $\Jtw=\sla_1^{i/4} (J_1 \ot J_2)$. Denote by
$\Jhtw$ the modular conjugation of the canonical left invariant weight
on the dual of $(\Mtw,\detw)$. Write $\Utw:=\Jhtw \Jtw$. Then, we know
that
\begin{align}
(1 \ot 1 \ot \Utw) (\Si V_1^* \Si)_{13} \; Z^*_{34} W_{2,24} Z_{34}
(1 \ot 1 \ot \Utw^*) & = (1 \ot 1 \ot \Utw) \Wtw (1 \ot 1 \ot \Utw^*)
\notag \\ &=
\Si \Vtw \Si = Z_{34}^* W^*_{1,13} Z_{34} (\Si V_2 \Si)_{24} \; . \label{eq.stom}
\end{align}
Taking slices on the first two legs, it follows that $\Utw (x \ot 1)
\Utw^* = Z^* (u_1 x u_1^* \ot 1) Z$ for all $x \in \Mh_1'$, where we use again the notation
$u_i = \Jh_i J_i$. Applying $\Ad (1 \ot 1 \ot \Utw^*)$ to Equation
\eqref{eq.stom}, realizing that $\Utw$ and $\Utw^*$ only differ up to
a scalar and again taking slices, we also conclude that $\Utw (1 \ot y)
\Utw^* = Z^* (1 \ot u_2 y u_2^*) Z$ for all $y \in \Mh_2'$. We
conclude that $\Ad \Utw = \Ad (Z^* (u_1 \ot u_2))$ on $\Mh_1' \ot
\Mh_2'$. As we know $\Jtw$, it follows that
$$\Jhtw z \Jhtw = Z^* (\Jh_1 \ot \Jh_1) z (\Jh_1 \ot \Jh_2) Z \; ,$$
for all $z \in \Mh_1' \ot \Mh_2'$. If next $z \in M_1 \ot M_2$, the
left hand side equals $\Rtw(z^*)$, but the right hand side as well
because $\tw^{-1}(R_1 \ot R_2)(z^*) = \Rtw(z^*)$. We conclude that
$\Jhtw$ and $Z^* (\Jh_1 \ot \Jh_2)$ only differ up to a scalar.

Because $\Jhtw^2=1$, we also get $(Z^* (\Jh_1 \ot
\Jh_2))^2=1$. Writing this out, it follows that
$$(J_1 \ot J_2) J \Jh = \Jh J (J_1 \ot J_2) \; .$$
Using this commutation relation, we arrive at
$$\Jtw Z^* (\Jh_1 \ot \Jh_2) = \sla_1^{i/2}
\bigl(\frac{\nu_1}{\nu_2}\bigr)^{i/4} Z^* (\Jh_1 \ot \Jh_2) \Jtw \;
.$$
On the other hand,
$$\Jtw \Jhtw = \nutw^{-i/4} \Jhtw \Jtw = \bigl(\frac{\nu_1}{\nu_2}\bigr)^{i/4}
\Jhtw \Jtw \; .$$
We conclude that $$\Jhtw = \sla_1^{i/4} Z^* (\Jh_1 \ot \Jh_2) =
\sla_1^{i/4} (J_1 \Jh_1 \ot J_2 \Jh_2) \Jh J (\Jh_1 \ot \Jh_2) =
\sla_1^{i/4} (\Jh_1 \ot \Jh_2) J \Jh (\Jh_1 J_1 \ot \Jh_2 J_2) \; .$$

Because we know $(\tautw_t)$ and the GNS-construction for $\vfitw$, we
immediately get that $\Ptw = P_1^{-1} \ot P_2$. Because
$$\nab_\tw^{it} = \nabp_1^{it} \ot \nab_2^{it} k_2^{it} J_2 k_2^{it}
J_2 \quad\text{and}\quad \nabp_\tw^{it} = \nab_1^{it} k_1^{it} J_1
k_1^{it} J_1 \ot \nabp_2^{it} \; ,$$
where $\nabp_1,\nabp_2$ and $\nabp_\tw$ denote the modular operators
of $\psi_1,\psi_2$ and $\psitw$ respectively, it is easy to find, out of
Proposition 2.4 in \cite{VVD}, the formula
$$\sdehtw^{it} = \Jh_1 \sdeh_1^{-it} \Jh_1 \; k_1^{-it} J_1 k_1^{-it} J_1 \ot \sdeh_2^{it} \; .$$
But, from Proposition \ref{prop.rho}, it follows that $\Jh (1 \ot \sdeh_2^{it}) \Jh = k_1^{-it} J_1 k_1^{-it} J_1 \ot \Jh_2 \sdeh_2^{it} \Jh_2$.
Because $1 \ot \sdeh_2^{it}$ is a group-like unitary of $(\Mh,\deh)$, we know that $J(1 \ot \sdeh_2^{it})J = 1 \ot \sdeh_2^{it}$ and so
$\sdehtw^{it} = (\Jh_1 \sdeh_1^{-it} \Jh_1 \ot 1) Z^*(1 \ot \sdeh_2^{it})Z$. The formula for $\Jhtw \sdehtw^{it}\Jhtw$ is proven analogously.
\end{proof}

\section{Actions of double crossed products and bicrossed products} \label{sec.actions}

Suppose that $\rho : N \recht \Mtw \ot N$ is an action of the double crossed product $(\Mtw,\detw)$ on the von Neumann algebra $N$. We will show how
to construct another von Neumann algebra $L$ with an action of the bicrossed product $(M,\de)$ and we show that the outerness of one action is
equivalent to the outerness of the other action. As a definition of outerness, we take Definition 5.5 in \cite{V}: $\Mtw \kruisje{\rho} N \cap
\rho(N)' = \C$. In \cite{V3}, the second author proves that every
l.c.\ quantum group can act outerly in this strong sense on a factor.

It is clear that we have morphisms $(\Mtw,\detw) \recht (M_1,\deoop)$
and $(\Mtw,\detw) \recht (M_2,\de_2)$. Although these morphisms do not
exist on the von Neumann algebra level, they can nevertheless be used
to {\it restrict} the action $\rho$ to an action of $(M_1,\deoop)$,
resp.\ $(M_2,\de_2)$. We briefly explain how to do this.

In the von Neumann algebra language, a morphism appears rather as a
special type of action (the {\it restricted} left or right
translation), as was pointed out by J. Kustermans in \cite{JK}. We
clearly have a left action $\zeta_1 : \Mtw \recht M_1 \ot \Mtw : \zeta_1
= \deoop \ot \io$ of $(M_1,\deoop)$ on $\Mtw$ satisfying $(\io \ot
\detw)\zeta_1 = (\zeta_1 \ot \io \ot \io)\detw$. Analogously, we have
the right action $\io \ot \de_2 : \Mtw \recht \Mtw \ot M_2$, which can
be composed with the unitary antipodes to yield the left action of
$(M_2,\de_2)$ on $\Mtw$ given by
$\zeta_2 : \Mtw \recht M_2 \ot \Mtw : \zeta_2 = \si_{12} \tw_{12} (\io
\ot \de_2)$.

This allows us to define the action $\eta_1 : N \recht M_1 \ot N$ of
$(M_1,\deoop)$ on $N$ such that $(\io \ot \rho)\eta_1 = (\zeta_1 \ot
\io)\rho$ and the action $\eta_2 : N \recht M_2 \ot N$ of $(M_2,\det)$
on $N$ such that $(\io \ot \rho)\eta_2 = (\zeta_2 \ot \io)\rho$.
Denote
$$L:= M_1 \kruisje{\eta_1} N = (\eta_1(N) \cup \Mh_1' \ot 1)\dpr \;
.$$
It is easy to check that $(\io \ot \eta_2)\eta_1 = \rho$ and applying
$\io \ot \eta_2$ to $L$, we arrive at another representation of $L$:
$(\rho(N) \cup \Mh_1' \ot 1 \ot 1)\dpr$. This shows how $L$ is an
intermediate von Neumann algebra
$$N \hookrightarrow L \hookrightarrow \Mtw \kruisje{\rho} N \; .$$
On $L$, we will define a natural action $\eta$ of the bicrossed product
$(M,\de)$ such that
\begin{equation} \label{eq.inclusion}
N \hookrightarrow L \hookrightarrow \Mtw \kruisje{\rho} N
\hookrightarrow M \kruisje{\eta} L \hookrightarrow \Mh_\tw
\kruisje{\hat{\rho}} (\Mtw \kruisje{\rho} N) \; .
\end{equation}

\begin{proposition}
There exists a unique action $\eta : L \recht M \ot L$ of $(M,\de)$ on
$L$ such that
$$\eta(\eta_1(x)) = (\al \ot \eta_1)\eta_2(x) \quad\text{for all}\;\;
x \in N \qquad \eta(\Jh_1 y \Jh_1 \ot 1) = (\Jh \ot \Jh_1)\mu(y) (\Jh
\ot \Jh_1) \ot 1 \quad\text{for all}\;\; y \in \Mh_1 \; ,$$
where $\mu : \Mh_1 \recht M \ot \Mh_1$ is the left action of $(M,\deop)$
on $\Mh_1$ defined by $\mu(y) \ot 1 = \deop(y \ot 1)$.

Moreover, the action $\eta$ is outer if and only if the action $\rho$ is outer. Finally, the inclusion of Equation~\eqref{eq.inclusion} holds.
\end{proposition}

\begin{proof}
During the proof, we replace $L$ by its representation $(\rho(N) \cup
\Mh_1' \ot 1 \ot 1)\dpr$ in $\B(H_1 \ot H_2) \ot N$. Define for $z \in L$,
$$\eta(z) = \cU_{123} (\io \ot \io \ot \rho)(z) \cU_{123}^*
\qquad\text{where}\quad \cU = (\Jh \ot \Jh_1) (\be \ot \io)(W_1) (\Jh_1 \ot \Jh_2 \ot \Jh_1) \; .$$ Represent $N$ on its standard Hilbert space, with
anti-unitary $J_N$ and use the notation $$\rhotil(y) = (\Jh_1 \ot \Jh_2 \ot J_N)\rho(J_N y J_N) (\Jh_1 \ot \Jh_2 \ot J_N) \quad\text{for}\quad y \in
N' \; .$$ Let $x \in N$. We have
\begin{align*}
& \eta(\rho(x)) = \cU_{123} (\detw \ot \io)\rho(x) \cU_{123}^* \\ &=
(\Jh \ot \Jh_1 \ot \Jh_2 \ot J_N) (\be \ot \io)(W_1)_{123} \; \bigl(\si_{23}
\tau_{23}^{-1} (\de_1 \ot \detop \ot \io) \rhotil(J_N x J_N) \bigr) \;
(\be \ot \io)(W_1^*)_{123} (\Jh \ot \Jh_1 \ot \Jh_2 \ot J_N) \\ &=
(\Jh \ot \Jh_1 \ot \Jh_2 \ot J_N) (\al \ot \io \ot \io \ot \io) \bigl(
\si_{12} \tau_{12}^{-1} (\io \ot \detop \ot \io)\rhotil(J_N x J_N)
\bigr) (\Jh \ot \Jh_1 \ot \Jh_2 \ot J_N) \\ &=
(\al \ot \io \ot \io \ot \io) \bigl( \si_{12} \tau_{12} (\io \ot \det
\ot \io)\rho(x) \bigr) = (\al \ot \rho)\eta_2(x) \; .
\end{align*}
On the other hand, it is obvious that $\eta(\Jh_1 y \Jh_1 \ot 1) =
(\Jh \ot \Jh_1)\mu(y) (\Jh \ot \Jh_1) \ot 1$ for $y \in \Mh_1$. Then,
it is also clear that $\eta : L \recht M \ot L$ and that $\eta$ is an
action of $(M,\de)$ on $L$.

By definition, $M \kruisje{\eta} L = (\eta(L) \cup \Mh \ot 1 \ot 1 \ot
1)\dpr$, but using $\cU$, we easily identify
$$M \kruisje{\eta} L = (L \cup (\Jh_1 \ot \Jh_2)\Jh \Mh \Jh (\Jh_1 \ot
\Jh_2) \ot 1)\dpr \subset \B(H_1 \ot H_2) \ot N \; .$$
In this way, we see that the inclusion of
Equation~\eqref{eq.inclusion} is given by
$$\rho(N) \subset L \subset (\rho(N) \cup \Mh_\tw \ot 1)\dpr \subset
(L \cup (\Jh_1 \ot \Jh_2)\Jh \Mh \Jh (\Jh_1 \ot \Jh_2) \ot 1)\dpr \subset
\B(H_1 \ot H_2) \ot N \; .$$

Suppose first that $\rho$ is outer and let $z \in (M \kruisje{\eta} L)
\cap L'$. Then, we observe that $z \in \B(H_1 \ot H_2) \ot N \cap
\rho(N)'$, which equals $M_1' \ot M_2' \ot 1$ by the outerness of
$\rho$. So, $z = y \ot 1$ with $y \in M_1' \ot M_2'$. But $z$ also
commutes with $\Mh_1' \ot 1 \ot 1$, which yields $z = 1 \ot a \ot 1$
with $a \in M_2'$. Because
$$(\io \ot \io \ot \eta_2)\rho = (\io \ot \io \ot \eta_2) (\io \ot
\eta_2) \eta_1 = (\io \ot \de_2 \ot \io)\rho \; ,$$
we see that $(\io \ot \io \ot \eta_2)(z) = V_{2,23} z_{124}
V_{2,23}^*$ for all $z \in M \kruisje{\eta} L$. So, $a \in \Mh_1$ and
hence, $a \in \C$. We get that $\eta$ is outer.

Suppose next that $\eta$ is outer and let $z \in \B(H_1 \ot H_2) \ot N
\cap (\Mtw \kruisje{\rho} N)'$. We claim that,
$$\cU_{123} (\io \ot \io \ot \rho)(y) \cU_{123}^* \in
\B(H_1 \ot H_2) \ot L \quad\text{for all}\quad y \in \B(H_1 \ot H_2) \ot N \; .$$ When $y \in \rho(N)$, this is clear and when $y \in \B(H_1 \ot H_2)
\ot 1$, it follows because $\cU \in \B(H_1 \ot H_2) \ot \Mh_1'$. Because $\rho$ is an action, we know that $$\B(H_1 \ot H_2) \ot N = (\rho(N) \cup
\B(H_1 \ot H_2) \ot 1)\dpr$$ and our claim is proven. So, we find in particular that $\cU_{123} (\io \ot \io \ot \rho)(z) \cU_{123}^* \in \B(H_1 \ot
H_2) \ot L \cap \eta(L)'$. From the outerness of $\eta$, it follows that $\B(H_1 \ot H_2) \ot L \cap \eta(L)' = M' \ot 1 \ot 1 \ot 1$ and we can take
$a \in M'$ such that
$$(\io \ot \io \ot \rho)(z) = \cU_{123}^* (a \ot 1 \ot 1 \ot 1)
\cU_{123} \; .$$
But, $\rho(N) \subset M_1 \ot M_2 \ot N$, while the last leg of $\cU$
sits in $\Mh_1'$. Hence, $z = y \ot 1$ for some $y \in \B(H_1 \ot H_2)$.
Writing $x = (\Jh_1 \ot \Jh_2)y(\Jh_1 \ot \Jh_2)$ and $b = \Jh a \Jh
\in M'$, we get
$$(\be \ot \io)(W_1) (x \ot 1) (\be \ot \io)(W_1^*) = b \ot 1 \; .$$
If we write $\zeta(c) = (\be \ot \io)(W_1) (c \ot 1) (\be \ot
\io)(W_1^*)$ for $c \in \B(H_1 \ot H_2)$, we see that $\zeta : \B(H_1
\ot H_2) \recht \B(H_1 \ot H_2) \ot \Mh_1$ is a right action of
$(\Mh_1,\dehoop)$ on $\B(H_1 \ot H_2)$. Because $\zeta(x) = b \ot 1$,
we find that $\zeta(b) \ot 1 = (\zeta \ot \io)\zeta(x) = (\io \ot \io
\ot \dehoop)\zeta(x) = b \ot 1 \ot 1$. Hence, $\zeta(b) = b \ot 1$ and
$x =b$. We arrive at the conclusion that $x \in M' \cap \be(M_1)' =
M_1' \ot M_2' \cap (\Mh_1 \ot 1)' = 1 \ot M_2'$. This gives us that $z
= 1 \ot d \ot 1$ with $d \in M_2'$. Because $z \in (\Mtw \kruisje{\rho} N)'$, $1 \ot d$ must commute with
$Z^*(1 \ot \Mh_2) Z$. Up to a scalar, $Z^*$ equals $(\Jh_1 J_1 \ot 1)
U_\be^* \; J(J_1 \ot J_2)$. The first part belongs to $\B(H_1) \ot M_2$
and we find that $1 \ot d$ must commute with $J(1 \ot \Mh_2)J$. This
means that $\al(J_2 d J_2)$ commutes with $1 \ot \Mh_2$. Hence,
$\al(J_2 d J_2) \in M_1 \ot 1$. Because $\al$ is an action of
$(M_1,\de_1)$ on $M_2$, it follows that $J_2 d J_2 \in \C$. So, $z \in
\C$ and we are done.
\end{proof}

\section{A characterization of double crossed products} \label{sec.char}

In Definition \ref{defin.matching}, we defined a matching $\tw$ of
l.c.\ quantum groups $(M_1,\de_1)$, $(M_2,\det)$ and we constructed the associated
double crossed product $(\Mtw,\detw)$ as a l.c.\ quantum group in
Section \ref{sec.doublelc}. It is clear that the mappings
$\deoop \ot \io$ and $\io \ot \det$ define, respectively, a left
action of $(M_1,\deoop)$ and a right action of $(M_2,\de_2)$ on the
von Neumann algebra $\Mtw \cong M_1 \ot M_2$. In fact, these actions
are precisely the actions that correspond to the natural morphisms
from $(\Mtw,\detw)$ to $(M_1,\deoop)$ and $(M_2,\det)$.

We claim that this characterizes
double crossed products: whenever we have a l.c.\ quantum group
$(N,\Gamma)$ and morphisms from $(N,\Gamma)$ to $(M_1,\deoop)$ and
$(M_2,\de_2)$, such that $N$, equipped with the naturally associated
left action of $(M_1,\deoop)$ and
right action of $(M_2,\de_2)$, is isomorphic to $M_1 \ot
M_2$ as an $M_1,M_2$-bicomodule, we get that $(N,\Gamma)$ is a double
crossed product.

We can interpret this in the classical, commutative case as
follows. We are given l.c.\ groups $G,G_1$ and
$G_2$. We suppose that $\Goop$ acts on the left and $G_2$ on the right
on $G$ such that, as a measure space, $G$ is isomorphic to $G_1 \times
G_2$ where the isomorphism intertwines the actions of $\Goop$ and
$G_2$ on $G$ with the natural actions of $\Goop$ and
$G_2$ on $G_1 \times G_2$. So, we assume in fact that $\Goop$ and
$G_2$ are closed subgroups of $G$ and that there exists an
element $x \in G$ such that the map $G_1 \times G_2 \rightarrow G :
(g,s) \mapsto gxs$ is a measure class isomorphism. Then, it is clear
that, considering $x^{-1} G_1 x$ and $G_2$ as closed subgroups of $G$,
we arrive at a matched pair (in the sense of \cite{BSV}, Definition
3.1) of l.c.\ groups with $G$ as the double crossed product. It is
obvious that it is really necessary to allow for the presence of $x
\in G$ and this explains the presence of the group-like unitary $u \in
\Nh$ in Proposition \ref{prop.characterization}.

\begin{proposition} \label{prop.characterization}
Let $(N,\Ga)$, $(M_1,\de_1)$ and $(M_2,\de_2)$ be l.c.\ quantum
groups. Suppose that there are morphisms from $(N,\Gamma)$ to
$(M_1,\deoop)$ and $(M_2,\det)$ respectively, with associated actions
$\zeta_1 : N \recht M_1 \ot N$, $\zeta_2 : N \recht N \ot M_2$ of
$(M_1,\deoop)$ on the left and $(M_2,\det)$ on the right on
$N$. Suppose that there exists an isomorphism of von Neumann algebras
$\pi : N \recht M_1 \ot M_2$ intertwining $\zeta_1$ with $\deoop \ot
\io$ and $\zeta_2$ with $\io \ot \det$.

Then, there exists a matching $\tw: M_1 \ot M_2 \recht M_1 \ot M_2$
and a group-like unitary $u \in \Nh$ such that $\pi \circ \Ad u : N
\recht M_1 \ot M_2$ is an isomorphism between the l.c.\ quantum groups
$(N,\Gamma)$ and the double crossed product $(\Mtw,\detw)$.
\end{proposition}
\begin{proof}
Denote by $W_N$ the left regular representation of $(N,\Gamma)$.
Because the actions $\zeta_1,\zeta_2$ are defined by morphisms, we
know that
\begin{align}
& (\zeta_1 \ot \io)\Gamma = (\io \ot \Gamma)\zeta_1 \; , \quad (\io \ot
\zeta_2)\Gamma = (\Gamma \ot \io)\zeta_2 \; , \notag \\
& (\zeta_1 \ot \io)(W_N) = X_{1,13} W_{N,23} \; , \quad (\zeta_2 \ot
\io)(W_N) = W_{N,13} X_{2,23} \; , \label{eq.zetas}
\end{align}
where $X_1 \in M_1 \ot \Nh$ is a corepresentation of $(M_1,\deoop)$ and
$X_2 \in M_2 \ot \Nh$ is a corepresentation of $(M_2,\det)$. Defining
$\Wtil_N = (\pi \ot \io)(W_N) \in M_1 \ot M_2 \ot \Nh$, we deduce that
$$(\deoop \ot \io \ot \io)(\Wtil_N) =X_{1,14} \Wtil_{N,234} \; , \quad
(\io \ot \det \ot \io)(\Wtil_N) = \Wtil_{N,124} X_{2,34} \; .$$
It then follows that $X^*_{1,13} \Wtil_N$ is invariant under $\deoop
\ot \io \ot \io$ and we find $Y \in M_2 \ot \Nh$ such that $\Wtil_N =
X_{1,13} Y_{23}$. Using next the formula for $(\io \ot \det \ot
\io)(\Wtil_N)$, we get $(\det \ot \io)(Y) = Y_{13} X_{2,23}$. Hence,
$YX^*_2$ is invariant under $\det \ot \io$ and we find a unitary $u
\in \Nh$ such that
$$\Wtil_N = X_{1,13} (1 \ot 1 \ot u) X_{2,23} \; .$$
We know that $(\io \ot \io \ot \deh_N)(\Wtil_N) = \Wtil_{N,124}
\Wtil_{N,123}$, from which we conclude that $(\io \ot \deh_N)(X_i) =
X_{i,13} X_{i,12}$ for $i=1,2$. Writing this out, it follows that
$\deh(u) = u \ot u$.

The equation $\deh(u) = u \ot u$ gives us that $(\Ad u \ot \io)(W_N) =
(1 \ot u^*)W_N$. From Equation \eqref{eq.zetas}, it is clear that
$\zeta_2$ commutes with the isomorphism $\Ad u$ of $N$, while defining
$\zetatil_1 = (\io \ot \Ad u^*) \circ \zeta_1 \circ \Ad u$ and $\Xtil_1 = (1
\ot u^*)X_1 (1 \ot u)$, we get
$$(\zetatil_1 \ot \io)(W_N) = \Xtil_{1,13} W_{N,23} \; .$$
Finally, we put $\pitil = \pi \circ \Ad u$ and observe that $(\pitil
\ot \io)(W_N)=\Xtil_{1,13} X_{2,23}$. Also, $\pitil$ intertwines
$\zetatil_1$ with $\deoop \ot \io$ and $\zeta_2$ with $\io \ot \det$.
As a conclusion, we see that
changing the morphism from $(N,\Gamma)$ to $(M_1,\deoop)$ by using
$u$, we may in fact suppose from the beginning that $u=1$ and now
prove the existence of a matching $\tw : M_1 \ot M_2 \recht M_1 \ot
M_2$ such that $\pi$ is an isomorphism between the l.c.\ quantum
groups $(N,\Gamma)$ and $(\Mtw,\detw)$.

So, we assume that $u=1$ and get $\Wtil_N = X_{1,13}
X_{2,23}$. Because a morphism commutes with the antipode, we know that
$\zeta_1 \tau^N_t = (\tau^1_{-t} \ot \tau^N_t) \zeta_1$ and $\zeta_2
\tau^N_t = (\tau^N_{t} \ot \tau^2_t) \zeta_2$. Combining this with
Equation \eqref{eq.zetas}, we get that $X_1$ is invariant under
$\tau^1_{-t} \ot \tauh^N_t$, while $X_2$ is invariant under $\tau^2_t
\ot \tauh^N_t$. Because $\Wtil_N = X_{1,13}
X_{2,23}$, it follows that $\pi \tau^N_t = (\tau^1_{-t} \ot
\tau^2_t)\pi$. Moreover, because $X_1$ and $X_2$ are
corepresentations, it follows that, for all $\om \in \Nh_*$, the
element $y:= (\io \ot \io \ot \om)(X_{1,13}X_{2,23}) \in D(S_1^{-1}
\ot S_2)$ and
$$(S_1^{-1} \ot S_2)(y) = (\io \ot \io \ot \om)(X^*_{1,13}X^*_{2,23})
\; .$$
Define $\tw := (R_1 \ot R_2) \pi R_N \pi^{-1}$, use that $S_N(\io \ot
\om)(W_N) = (\io \ot \om)(W_N^*)$ and combine this with the formulas
above, to conclude that
$$(\tw \ot \io)(X_{1,13}X_{2,23}) = X_{2,23}X_{1,13} \; .$$
It follows that $\bigl((\io \ot \si\tw \ot \io)(\deoop \ot \det) \ot
\io\bigr)(\Wtil_N) = \Wtil_{N,125} \Wtil_{N,345}$ and hence, $(\io \ot
\si\tw \ot \io)(\deoop \ot \det) \pi = (\pi \ot \pi)\Gamma$. We also
verify easily that
\begin{align*}
& ((\de_1 \ot \io)\tw \ot \io)(X_{1,13} X_{2,23}) = X_{2,34} X_{1,24}
X_{1,14} \\
\text{and} \quad & \tw_{23} \tw_{13} (\de_1 \ot \io \ot \io)(X_{1,13}
X_{2,23}) = \tw_{23} \tw_{13} (X_{1,24} X_{1,14} X_{2,34}) = X_{2,34} X_{1,24}
X_{1,14} \; .
\end{align*}
Because the slices of $\Wtil_N$ are dense in $\pi(N) = M_1 \ot M_2$,
we conclude that $(\de_1 \ot \io)\tw = \tw_{23} \tw_{13} (\de_1 \ot
\io)$. Completely analogously, we verify that $(\io \ot \de_2)\tw =
\tw_{13} \tw_{12} (\io \ot \de_2)$. Hence, $\tw$ is a matching and
above, we have seen already that $\pi : N \recht \Mtw$ is an
isomorphism of l.c.\ quantum groups.
\end{proof}

\section{Example: generalized quantum doubles} \label{sec.doublegen}

Let $(\Mo,\deo)$ and $(\Mt,\det)$ be l.c.\ quantum groups. Suppose
that $\cZ \in \Mo \ot \Mt$ is a bicharacter in the sense that
$$(\deo \ot \io)(\cZ) = \cZ_{23} \cZ_{13} \quad\text{and}\quad (\io
\ot \det)(\cZ) = \cZ_{13} \cZ_{12} \; .$$
(We see that, in fact, $\cZ^*$ is the bicharacter.) Then, we can
define a matching $\tw$ as the inner isomorphism:
\begin{equation} \label{eq.inner}
\tw(x) = \cZ x \cZ^* \quad\text{for all}\quad x \in \Mo \ot \Mt \;
.
\end{equation}
It is clear that $\tw$ is a matching of $(\Mo,\deo)$ and
$(\Mt,\det)$. Hence, we have an associated bicrossed product $(M,\de)$
and double crossed product $(\Mtw,\detw)$, as above.

A particular case is that of the \emph{quantum double}: let
$(N,\de_N)$ be a l.c.\ quantum group and put $(M_1,\deo) =
(N,\deopN)$, $(M_2,\det) = (\Nh,\deh_N)$ and $\cZ = W_N$. The
associated double crossed products is the quantum double of
$(N,\de_N)$. Such a quantum double construction in the framework of
l.c.\ quantum groups (or rather, Woronowicz algebras) has been done
for the first time by T. Yamanouchi \cite{Yam3}.

For this reason, we will call the double crossed product associated
with an inner matching as in Equation \eqref{eq.inner} a
\emph{generalized quantum double}.

Suppose for the rest of this section that we have a bicharacter $\cZ$
as above and an inner matching $\tw$ defined by Equation
\eqref{eq.inner}.

\begin{proposition}
Through the isomorphism $\Ad \cZ^*$, the bicrossed product $(M,\de)$
is isomorphic with the tensor product $(\Mh_1 \ot \Mt, (\io \ot \si
\ot \io)(\deh_1 \ot \det))$.

The ingredients of the double crossed product (i.e., the generalized
quantum double) are given by Theorem \ref{thm.double}, combined with
the formulas:
\begin{itemize}
\item $J = \cZ (\Jh_1 \ot J_2) \cZ^* \quad\text{and}\quad \Jh = \cZ
  (J_1 \ot \Jh_2) \cZ^*$;
\item there exist group-like unitaries $k_1^{it} \in M_1$ and
  $k_2^{it} \in M_2$ such that $$\cZ^*(\sdeh_1^{it} \ot 1) \cZ =
  \sdeh_1^{it} \ot k_2^{-it} \quad\text{and}\quad \cZ^*(1 \ot
  \sdeh_2^{it}) \cZ = k_1^{-it} \ot \sdeh_2^{it}$$ and $\rho_1^{it} =
  k_1^{it} \sde_1^{it}$, $\rho_2^{it} = k_2^{it} \sde_2^{it}$ (observe
  that $\rho_1$ and $\rho_2$ are defined in Proposition
  \ref{prop.rho}, associated with a matching $\tw$);
\item $Z = \nu_1^{i/4} \cZ \; (J_1 \ot J_2)\cZ(J_1 \ot J_2)$.
\end{itemize}

The double crossed product is unimodular if and only if $k_1=\sde_1$
and $k_2=\sde_2$. This is the case for the ordinary quantum double of
a l.c.\ quantum group.
\end{proposition}
\begin{proof}
To determine the bicrossed product, we calculate:
\begin{align*}
(\cZ^* \ot \cZ^*)W (\cZ \ot \cZ) &= \cZ^*_{34} W_{2,24} \cZ^*_{12}
\cZ_{34} \Wh_{1,13} \cZ_{12} \\ &= W_{2,24} (\io \ot \det)(\cZ^*)_{324}
\cZ^*_{12} \cZ_{34} (\deoop \ot \io)(\cZ)_{132} \Wh_{1,13} = W_{2,24}
\Wh_{1,13} \; .
\end{align*}
This proves that $\Ad \cZ^*$ maps $(M,\de)$ onto the tensor product
$(\Mh_1 \ot M_2, (\io \ot \si \ot \io)(\deh_1 \ot \det))$. From
Proposition 4.2 in \cite{V}, the formulas for $J$ and $\Jh$ follow.

Because $\cZ^*$ is a corepresentation for $(M_1,\de_1)$, we get that
$(\io \ot \om)(\cZ^*) \in D(S_1)$ for all $\om \in M_{2*}$ and $S_1 (\io
\ot \om)(\cZ^*) = (\io \ot \om)(\cZ)$. An analogous statement holds
for the left slices of $\cZ$ and $S_2$. Combining both, we conclude
that $\cZ$ is invariant under $\tauo_{-t} \ot \taut_t$ and $R_1 \ot
R_2$. The formula for $Z$ then follows.

Finally, $\etao_t(x) = \sdeh_1^{it} x \sdeh_1^{-it}$ defines a
one-parameter group of automorphisms of $M_1$ satisfying $\deo \etao_t
= (\etao_t \ot \io)\deo$. It follows that $(\etao_t \ot \io)(\cZ^*)
\cZ$ is invariant under $\deo \ot \io$. Hence, there exist unitaries
$k_2^{it} \in M_2$ such that $\cZ^*(\sdeh_1^{it} \ot 1) \cZ =
\sdeh_1^{it} \ot k_2^{-it}$. Applying $\io \ot \det$ to this equality,
we find that $k_2^{it}$ is group-like. We analogously find
$k_1^{it}$. Because the modular element of the tensor product $\Mh_1
\ot M_2$ is given by $\sdeh_1 \ot \sde_2$, we get the following
formula for the modular element $\sde$ of $(M,\de)$:
$$\sde^{it} = \cZ (\sdeh_1^{it} \ot \sde_2^{it}) \cZ^* = (\sdeh_1^{it}
\ot 1) \al(k_2^{it} \sde_2^{it}) \; .$$ So, by definition,
$\rho_2^{it} = k_2^{it} \sde_2^{it}$.

It follows immediately from Theorem \ref{thm.double} that
$(\Mtw,\detw)$ is unimodular if and only if $k_1 = \sde_1$ and $k_2 =
\sde_2$. This is the case for the ordinary quantum double: we have
$(M_1,\deo) = (N,\deopN)$ and hence, $\sdeh_1^{it}= \Jh_N
\sdeh_N^{-it} \Jh_N$ satisfying
$$\cZ^* (\sdeh_1^{it} \ot 1)\cZ = (\Jh_N \ot J_N)W_N (\sdeh_N^{-it}
\ot 1) W_N^* (\Jh_N \ot J_N) = \Jh_N \sdeh_N^{-it} \Jh_N \ot
\sdeh_N^{-it} = \sdeh_1^{it} \ot \sde_2^{-it} \; .$$
So, we are done.
\end{proof}

\section{C$^\mathbf{*}$-algebraic features of bicrossed and double
  crossed products} \label{sec.cstar}

Throughout this section, we fix two l.c.\ quantum groups $(\Mo,\deo)$, $(\Mt,\det)$ and a matching $\tw$ on them. We have an associated bicrossed
product $(M,\de)$ and double crossed product $(\Mtw,\detw)$. We continue to use the notations fixed in Notation \ref{nota.crucial}. So, the
\cst-algebras associated to $(\Mo,\deo)$, $(\Mt,\det)$, $(M,\de)$ and $(\Mtw,\detw)$ are denoted by $A_1,A_2,A$ and $\Atw$, respectively.

We denote by $\K$ the \cst-algebra of compact operators on $H_1 \ot
H_2$ and by $\K_i$ the \cst-algebra of compact operators on $H_i$, $i=1,2$.

The following result is rather easy to prove. We give formulas for the \cst-algebras $A$, $\Ah$ associated with the bicrossed product and for the
\cst-algebra $\Ahtw$ associated with the double crossed product. Observe that it is much more delicate to describe the \cst-algebra $\Atw$, see
Proposition \ref{prop.atw} where we give a formula in the case of regular $(M_1,\de_1)$ and $(M_2,\de_2)$.

\begin{proposition} \label{prop.cstar}
We have
$$A=[(\Ah_1 \ot 1) \al(A_2)] \quad\text{and}\quad \Ah = [(1 \ot \Ah_2)
\be(A_1) ] \; .$$ Further, we have
$$\Ahtw=[ (\Jh_1 \Ah_1 \Jh_1 \ot 1) \; Z^* (1 \ot \Ah_2)Z ] \; .$$
\end{proposition}
\begin{proof}
We have
\begin{align*}
A &= [(\io \ot \io \ot \om \ot \mu)\bigl( (\al \ot \io \ot
\io)(W_{2,13} \; U_{\be,23}) \; \Wh_{1,13} \bigr)] \\ &= [(\io \ot \io \ot
\om \ot \mu)\bigl( (\al \ot \io \ot \io)( U_{\be,23} \; W_{2,13} \;
U_{\be,21}^* ) \; \Wh_{1,13} \bigr) ] \\ &= [ \al(A_2) (\io \ot \io \ot
\om)\bigl( (\al \ot \io)(U_{\be,21}^*) \; \Wh_{1,13} \bigr) ] =
[\al(A_2) (\Ah_1 \ot 1) ] \; ,
\end{align*}
where we used that $(\io \ot \det)(U_\be) = U_{\be,13} U_{\be,12}$ and
the fact that $U_\be \in \M(\K_1 \ot A_2)$.

We know that $A$ is a \cst-algebra. This proves the formula for $A$ and the formula for $\Ah$ is proven analogously. The formula for $\Ahtw$ is an
immediate corollary of the formula for $\Wtw$ in Theorem \ref{thm.double}.
\end{proof}

Next, we present a quantum version of Theorem 3.11 in \cite{BSV}: we
give a necessary and sufficient condition for $(M,\de)$ to be
(semi-)regular.

\begin{proposition}
We have that $(M,\de)$ is semi-regular if and only if both
$(M_1,\deo)$ and $(\Mt,\det)$ are semi-regular and $A_1 \ot A_2
\subset \Atw$.

We have that $(M,\de)$ is regular if and only if both $(M_1,\deo)$ and $(\Mt,\det)$ are regular and one of the following conditions holds true: $A_1
\ot A_2= \Atw$ or $\tw(A_1 \ot A_2) = A_1 \ot A_2$.
\end{proposition}
\begin{proof}
First, observe that
$$[\cC(W)] = [(\io \ot \io \ot \om \ot \mu)(\Si_{13} \Si_{24} W_{2,24}
U_{\be,34} U_{\al,12}^* \Wh_{1,13} )] \; .$$
Because $(\io \ot \det)(U_\be) = U_{\be,13} U_{\be,12}$ and $(\deo \ot
\io)(U_\al) = U_{\al, 23} U_{\al,13}$, we get
$$[\cC(W)] = U_\be \; [(\io \ot \io \ot \om \ot \mu)\bigl( (\Si
W_2)_{24} \; (U_\be^* U_\al)_{12} \; (\Si \Wh_1)_{13} \bigr)] \;
U^*_\al = U_\be \; [(1 \ot \cC(W_2)) \; U_\be^* U_\al \; (\cC(\Wh_1)
\ot 1)] \; U^*_\al \; .$$
Suppose that $(M,\de)$ is semi-regular. From Proposition 5.7 in
\cite{BSV}, it follows that semi-regularity passes to closed quantum
subgroups. So, $(M_1,\deo)$ and $(M_2,\det)$ are semi-regular. From
the calculation above, we conclude that
$$\K \subset [(1 \ot \cC(W_2)) \; U_\be^* U_\al \; (\cC(\Wh_1)
\ot 1)]$$
and hence,
$$\K = [(1 \ot \K_2)\K(\K_1 \ot 1)] \subset [(1 \ot \K_2) \; U_\be^* U_\al \; (\K_1
\ot 1)] \; .$$
Hence, $\K \subset [(1 \ot \K_2)Z^*(\K_1 \ot 1)]$.

On the other hand, we have
\begin{align*}
\Atw &= [(\io \ot \io \ot \om \ot \mu)\bigl( (\Si V_1^* \Si)_{13} \;
Z^*_{34} \; W_{2,24} \bigr)] = [(\io \ot \io \ot \om \ot \mu)\bigl(
(\Si V_1^* \Si)_{13} \; ((1 \ot \K_2)Z^*(\K_1 \ot 1))_{34} \; W_{2,24}
\bigr)] \\ &\supset [(\io \ot \io \ot \om \ot \mu)\bigl(
(\Si V_1^* \Si)_{13} \; (\K_1 \ot \K_2)_{34} \; W_{2,24}
\bigr)] = A_1 \ot A_2 \; .
\end{align*}

Conversely, suppose that $(M_1,\deo)$ and $(M_2,\det)$ are both
semi-regular and $A_1 \ot A_2 \subset \Atw$. It is easy to check that
$A_1 \ot A_2 = [\Atw (A_1 \ot 1)]$. Hence,
$$\tw(A_1 \ot A_2) = [(R_1 \ot R_2)\Rtw (\Atw) \; \be(A_1) ] \supset
[(A_1 \ot A_2)\be(A_1)] \; .$$
But, using Proposition \ref{prop.cstar}, we have
$$[A \; \Ah] = [(\Ah_1 \ot 1) \tw(A_1 \ot A_2) (1 \ot \Ah_2) ] \supset
[(\Ah_1 A_1 \ot A_2) \be(A_1) (1 \ot \Ah_2)] \supset [(\K_1 \ot A_2)
\; \be(A_1) \; (1 \ot \Ah_2)] \; .$$
Because $[A \; \Ah]$ is a \cst-algebra, we get
$$[A \; \Ah] \supset \K_1 \ot [A_2 \; (\om \ot \io)\be(A_1) \; \Ah_2]
= \K_1 \ot [A_2 \; \Ah_2] \supset \K_1 \ot \K_2 \; ,$$
where we used that $\be(x) = U_\be (x \ot 1) U_\be^*$ with $U_\be \in
\M(\K_1 \ot A_2)$. Hence, $(M,\de)$ is semi-regular. This concludes
the proof of the first statement of the proposition.

In order to prove the second statement, we can replace above every inclusion by an equality to obtain the equivalence between regularity of $(M,\de)$
and regularity of $(M_1,\de_1)$, $(M_2,\de_2)$ and the equality $\Atw=A_1 \ot A_2$. The only different point is to prove that the regularity of
$(M,\de)$ implies the regularity of $(M_1,\deo)$ and $(M_2,\det)$. But, suppose that $(M,\de)$ is regular. Then, $(\al \ot \io)(W_1)$ is a
corepresentation of $(M,\de)$ and satisfies as such
$$[(\K_1 \ot \K_2 \ot 1) (\al \ot \io)(W_1) (1 \ot 1 \ot \K_1)] = \K_1
\ot \K_2 \ot \K_1 \; .$$
Because $\al(x)= U_\al(1 \ot x)U_\al^*$, it follows immediately that
$(M_1,\deo)$ is regular. We get the regularity of $(M_2,\det)$ in an
analogous way.

Finally, we suppose that $(M_1,\de_1)$ and $(M_2,\de_2)$ are regular and $\tw(A_1 \ot A_2) = A_1 \ot A_2$. Using Proposition \ref{prop.cstar}, we
immediately have
$$[A \; \Ah] = [(\Ah_1 \ot 1) \tw(A_1 \ot A_2) (1 \ot \Ah_2) ] = [ \Ah_1 A_1 \ot A_2 \Ah_2 ] = \K_1 \ot \K_2 \; ,$$
yielding the regularity of $(M,\de)$.
\end{proof}

The (semi-)regularity of the double crossed product $(\Mtw,\detw)$ is
characterized more easily in the following proposition.

\begin{proposition}
We have that $(\Mtw,\detw)$ is semi-regular (resp., regular) if and
only if both $(M_1,\deo)$ and $(M_2,\det)$ are semi-regular (resp.,
regular).
\end{proposition}
\begin{proof}
First observe that
$$[\cC(\Wtw)] = [(\io \ot \io \ot \om \ot \mu)\bigl( (\Si V_1)^*_{13}
\; Z^*_{32} \; (\Si W_2)_{24} \bigr) ] = [(\io \ot \io \ot \om)\bigl((\Si V_1)^*_{13}
\; Z^*_{32}\bigr) \; (1 \ot \cC(W_2))] \; .$$
Suppose first that $(M_1,\deo)$ and $(M_2,\det)$ are
semi-regular. Then, it follows that
$$[\cC(\Wtw)] \supset [(\io \ot \io \ot \om)\bigl((\Si V_1)^*_{13}
\; Z^*_{32}\bigr) \; (1 \ot \K_2)] = [(\io \ot \io \ot \om)\bigl((\Si V_1)^*_{13}
\bigr) \; (1 \ot \K_2)] \supset \K_1 \ot \K_2 \; .$$
Hence, $(\Mtw,\detw)$ is semi-regular.

Suppose, conversely, that $(\Mtw,\detw)$ is semi-regular. Because
$(\Mh_1',\deh_1')$ and $(\Mh_2,\deh_2)$ are closed quantum subgroups
of $(\Mhtw,\dehtw)$, we get that they are semi-regular, using Proposition
5.7 in \cite{BSV}. Hence, the same holds for $(M_1,\deo)$ and
$(\Mt,\det)$.

To prove the statement about regularity, we proceed exactly as in the
final paragraph of the proof of the previous proposition.
\end{proof}

Finally, we want to give some more detailed information about the
\cst-algebras $\Atw$ and $A$, as well as the crossed product $A
\rtimes\red A$ of the right action of $(A,\de)$ on itself. Our results
will be quantum versions of the following results in \cite{BSV}:
Proposition 3.6, Lemma 3.9 and Proposition 3.10.

To obtain our results, we need actions of \cst-algebraic quantum
groups on \cst-algebras. As shown in Section 5 of \cite{BSV}, the
notion of continuity of such an action is problematic when the acting
quantum group is non-regular. So, for the rest of this section,
\emph{we suppose that $(M_1,\deo)$ and $(M_2,\det)$ are regular}.

\begin{proposition}
Defining
$$T_1 = [(\io \ot \mu)\be(A_1) \mid \mu \in \B(H_2)_*]
\quad\text{and}\quad T_2=[(\om \ot \io)\al(A_2) \mid \om \in \B(H_1)_*]
\; ,$$
we have that $T_1$ and $T_2$ are \cst-algebras, such that the
restrictions of $\al$ and $\be$:
$$\al : T_2 \recht \M(A_1 \ot T_2) \; , \quad \be : T_1 \recht \M(T_1
\ot A_2)$$
define \cst-algebraic actions of $(A_1,\deo)$ on the left on $T_2$ and
of $(A_2,\detop)$ on the right on $T_1$, which are continuous in the
strong sense:
$$[\al(T_2)(A_1 \ot 1)] = A_1 \ot T_2 \quad\text{and}\quad [\be(T_1)
(1 \ot A_2)] = T_1 \ot A_2 \; .$$
Moreover, we have
$$A= A_1 \ltimesred T_2 \quad\text{and}\quad \Ah = T_1 \rtimes\red
A_2 \; .$$
\end{proposition}
\begin{proof}
As in the proof of Proposition 5.7 in \cite{BSV}, we conclude that
$T_1$ and $T_2$ are \cst-algebras. Further, we observe that
\begin{align*}
[\al(T_2)(A_1 \ot 1)] &= [(\om \ot \io \ot \io)\bigl( V_{1,12}
\al(A_2)_{13} V_{1,12}^* \; (1 \ot A_1 \ot 1) \bigr) ] \\ &= [(\om \ot \io
\ot \io)\bigl( ((\K_1 \ot 1)V_1(1 \ot A_1))_{12} \; \al(A_2)_{13}
\bigr)] = A_1 \ot T_2 \; ,
\end{align*}
where we used that $V_1 = \M(\K_1 \ot A_1)$ and $[(\K_1 \ot 1)V_1(1
\ot A_1)] = \K_1 \ot A_1$ by regularity of $(M_1,\deo)$. It follows
that $\al : T_2 \recht \M(A_1 \ot T_2)$ is non-degenerate and that the restriction of
$\al$ to $T_2$ is continuous in the strong sense. We analogously show the same
kind of properties for $T_1$ and $\be$.

By definition, we have $A_1 \ltimesred T_2 = [(\Ah_1 \ot 1)\al(T_2)]$
and this a \cst-algebra. So,
\begin{align*}
A_1 \ltimesred T_2 &= [(\Ah_1 \ot 1) \al(T_2) (\Ah_1 \ot 1)] = [(\om
\ot \io \ot \io)\bigl( (1 \ot \Ah_1 \ot 1) \; W^*_{1,12} \al(A_2)_{23}
W_{1,12} \; (1 \ot \Ah_1 \ot 1) \bigr)] \\ &= [(\Ah_1 \ot 1) \al(A_2)
(\Ah_1 \ot 1)] \; ,
\end{align*}
because $W_1 \in \M(\K_1 \ot \Ah_1)$. But, from Proposition
\ref{prop.cstar}, we know that $A = [(\Ah_1 \ot 1)\al(A_2)]$ and we
know that $A$ is a \cst-algebra. Hence, we are done.
\end{proof}

Recall that comparing our work with \cite{BSV}, the double crossed
product $(\Mtw,\detw)$ should be considered as the \lq big\rq\ group
$G$ in \cite{BSV}. Hence, the following result is really the quantum
version of \cite{BSV}, Lemma 3.9 and Proposition 3.10.

\begin{proposition} \label{prop.atw}
We have
$$\Atw = [(\io \ot \io \ot \om)\bigl( \tw_{13} (\io \ot \detop)(A_1
\ot A_2) \bigr) \mid \om \in \B(H_2)_*] \; .$$
Further,
$$\deoop \ot \det : \Atw \recht \M(A_1 \ot \Atw \ot A_2)$$
defines a \cst-algebraic action of $(A_1,\deoop)$ on the left and of
$(A_2,\det)$ on the right on $\Atw$. This action is continuous in the
strong sense:
$$[(A_1 \ot 1 \ot 1 \ot 1) (\deoop \ot \det)(\Atw) (1 \ot 1 \ot 1 \ot
A_2)] = A_1 \ot \Atw \ot A_2 \; .$$
Moreover, the crossed product satisfies:
$$A \rtimes\red A \cong A_1 \ltimesred \Atw \rtimes\red A_2 \; .$$
\end{proposition}

\begin{proof}
One checks easily that $A_1 \ot A_2 = [\Atw (A_1 \ot 1)]$. Hence, it
follows that
\begin{align*}
[(\io \ot \io \ot \om)\bigl( \tw_{13}(\io \ot \detop)(A_1 \ot A_2)
\bigr)] &= [(\io \ot \io \ot \om)\bigl( \tw_{13}(\io \ot \detop)(\Atw)
\; \tw(A_1 \ot 1)_{13} \bigr)] \\ &= [(\io \ot \io \ot \om)\bigl(
\tw^{-1}_{12} ((1 \ot 1 \ot A_2)(\io \ot \detop)\tw(\Atw)) \; \tw(A_1
\ot 1)_{13} \bigr) ] \; .
\end{align*}
From Proposition \ref{prop.multun}, we get that
$$[(1 \ot 1 \ot A_2) (\io \ot \detop)\tw(\Atw)] = [(\io \ot \io \ot
\io \ot \om \ot \mu) \bigl( ((A_2 \ot 1)W_2)_{35} \; W_{2,25} \;
Z_{45} \; (\Si V_1^* \Si)_{14} \bigr) ] = \tw(\Atw) \ot A_2 \; ,$$
because $W_2  \in \M(A_2 \ot \K_2)$. Hence, we conclude that
\begin{equation}\label{eq.ster}
[(\io \ot \io \ot \om)\bigl( \tw_{13}(\io \ot \detop)(A_1 \ot A_2)
\bigr)] = [\Atw (T_1 \ot 1)] \; .
\end{equation}
Now, we use the second formula for $\Vtw$ in Theorem \ref{thm.double}
and we find that
\begin{equation} \label{eq.stertw}
[(T_1 \ot 1) \Atw] = [(\om \ot \io \ot \io)\bigl( ((T_1 \ot
\K_2)\be(A_1))_{21} \; U_{\al,21} \; V_{2,13} \bigr) ] \; .
\end{equation}
Next, we claim that $[(T_1 \ot \K_2)\be(A_1)] = [(1 \ot \K_2)
\be(A_1)]$. Because $\be$ is continuous in the strong sense on $T_1$,
we have $[(1 \ot \K_2)\be(T_1)] = T_1 \ot \K_2$. To prove our claim,
it suffices to show that $[T_1 \; A_1] = A_1$. This last statement can
be proven as follows: from Theorem \ref{thm.double}, we know that
$(\io \ot \tw)(U_{\be,13} \; V_{1,12}) = V_{1,12} \; U_{\be,13}$,
which gives
$$T_1 = [(\om \ot \io \ot \mu)( (\io \ot \tw)(U_{\be,13}^*) \;
V_{1,12}) ]$$
and so, because $V_1 \in \M(\K_1 \ot A_1)$ and $U_\al \in \M(A_1 \ot \K_2)$,
$$[T_1 \; A_1] = [(\om \ot \io \ot \mu)(\io \ot \al)(U_\be^*) \; A_1]
= [(\om \ot \io \ot \mu)(U_{\al,23} U_{\be,13}^* U_{\al,23}^*) \; A_1]
= A_1 \; .$$
Hence, our claim is proven. Combining this with Equation
\eqref{eq.stertw}, we get
$$[(T_1 \ot 1)\Atw] = [(\om \ot \io \ot \io)\bigl( ((1 \ot
\K_2)\be(A_1))_{21} \; U_{\al,21} \; V_{2,13} \bigr) ] = \Atw \; ,$$
where we used again the second formula for $\Vtw$ in Theorem
\ref{thm.double}. Combining this equation with Equation
\eqref{eq.ster}, we finally arrive at
$$[(\io \ot \io \ot \om)\bigl( \tw_{13}(\io \ot \detop)(A_1 \ot A_2)
\bigr)] = \Atw \; .$$

It is easy to check that $\deoop \ot \det$ defines a \cst-algebraic
action of $(A_1,\deoop)$ on the left and $(A_2,\det)$ on the right on
$\Atw$, which is continuous in the strong sense.

Because $\Rtw(\Atw) = \Atw$, we also get
$$\Atw = [(\io \ot \io \ot \om)(\io \ot \det)\tw^{-1}(A_1 \ot A_2)] \;
.$$
Now, by definition, $A \rtimes\red A = [A \; \Jh \Ah \Jh]$. Because
$\Jh A \Jh = A$, we get, using $\Ad ((\Jh_1 \ot \Jh_2)\Jh)$,
$$A \rtimes\red A \cong (\Jh_1 \ot \Jh_2) [A \; \Ah] (\Jh_1 \ot
\Jh_2)\; .$$
On the other hand, by definition,
$$A_1 \ltimesred \Atw \rtimes\red A_2 = [(\Jh_1 \Ah_1 \Jh_1 \ot 1 \ot
1 \ot \Jh_2 \Ah_2 \Jh_2) (\deoop \ot \det)(\Atw) ] \; ,$$
which, by using the unitary $\Wh_1^* \ot W_2$, is isomorphic with $[
(\Jh_1 \Ah_1 \Jh_1 \ot \Jh_2 \Ah_2 \Jh_2) \Atw]$.

Now, we observe that
$$(\Jh_1 \ot \Jh_2) [A \; \Ah] (\Jh_1 \ot \Jh_2) = [(\Jh_1 \Ah_1 \Jh_1
\ot 1) \; \tw^{-1}(A_1 \ot A_2) \; (1 \ot \Jh_2 \Ah_2 \Jh_2)] \; .$$
Because the left hand side is a \cst-algebra, we have
\begin{align*}
(\Jh_1 \ot \Jh_2) [A \; \Ah] (\Jh_1 \ot \Jh_2) &= [(\Jh_1 \Ah_1 \Jh_1
\ot \Jh_2 \Ah_2 \Jh_2) \; \tw^{-1}(A_1 \ot A_2) \; (1 \ot \Jh_2 \Ah_2
\Jh_2)] \\ &= [ (\Jh_1 \Ah_1 \Jh_1 \ot \Jh_2 \Ah_2 \Jh_2) \; (\io \ot
\io \ot \om)\bigl( \tw^{-1}(A_1 \ot A_2)_{12} \; V_{2,23}^* \bigr) ]
\\ &= [ (\Jh_1 \Ah_1 \Jh_1 \ot \Jh_2 \Ah_2 \Jh_2) \; (\io \ot
\io \ot \om)\bigl( V_{2,23}^* \; (\io \ot \det)\tw^{-1}(A_1 \ot A_2)
\bigr) ] \\ &= [ (\Jh_1 \Ah_1 \Jh_1 \ot \Jh_2 \Ah_2 \Jh_2) \; (\io \ot
\io \ot \om)(\io \ot \det)\tw^{-1}(A_1 \ot A_2) ] \\ &= [ (\Jh_1 \Ah_1
\Jh_1 \ot \Jh_2 \Ah_2 \Jh_2) \; \Atw] \; ,
\end{align*}
where we used twice that $V_2 \in \M(\Jh_2 \Ah_2 \Jh_2 \ot \K_2)$. We
have seen above that this last subspace is isomorphic with $A_1
\ltimesred \Atw \rtimes\red A_2$. So, we are done.
\end{proof}

\section{Appendix: the Radon-Nikodym derivative of a weight under a quantum
group action} \label{sec.app}

In this appendix, we generalize results of T.~Yamanouchi
\cite{Yam1,Yam2} and we give considerably simpler proofs.

Let $N$ be a von Neumann algebra with \nsf weight $\te$. If a l.c.\
group $G$ acts continuously by automorphisms $(\al_g)$ on $N$, we can
consider the Connes cocycles (Radon-Nikodym derivatives) $[D \te \na
\al_g : D \te]_t$ as a function
$\R \recht L^\infty(G) \ot N$. In the case of an action of a l.c.\
quantum group $(M,\de)$ on $N$, Yamanouchi \cite{Yam1} defined the analogous
Radon-Nikodym derivative $[D \te \na \al : D \te]_t \in M \ot N$. We
repeat his definition.

We first introduce some needed preliminaries.

Fix an action $\al : N
\recht M \ot N$ of the l.c.\ quantum group $(M,\de)$ on the von
Neumann algebra $N$ and fix an \nsf weight $\te$ on $N$.
Denote by $\Jtil_\te$ and $\nabtil_\te$ the modular operators of the dual
weight $\tetil$, acting on the Hilbert space $H \ot H_\te$.

On the crossed product $M \kruisje{\al} N$, we have a dual action of
$(\Mh,\dehop)$ and the second crossed product $\Mh \kruisje{\alh} (M
\kruisje{\al} N)$ is canonically isomorphic with $B(H) \ot N$ through
the isomorphism $\Phi : B(H) \ot N \recht \Mh \kruisje{\alh} (M
\kruisje{\al} N)$, see \cite{V}, Theorem~2.6. So, we can consider the
bidual weight $\tetiltil$ on $B(H) \ot N$, with canonical GNS-map
$$(xy \ot 1) \al(z) \mapsto \la'(x) \ot \lah(y) \ot \la_\te(z)
\quad\text{for}\quad x \in M', y \in \Mh, z \in N$$
and where $\la'(x) = J \la(JxJ)$ is the canonical GNS-map for the
left invariant weight on the commutant quantum group $(M',\de')$,
which is the dual of $(\Mh,\dehop)$. In the same sense as above, these
elements $(xy \ot 1) \al(z)$ span a core for the GNS-map of
$\tetiltil$.

From the results of \cite{VVD}, Section~2, it then follows that
$$\tetiltil(\al(a^*) (x^*y \ot 1) \al(b)) = \Tr_{\nabh}(x^*y)
\te(a^*b)$$
for all $x,y \in \cN_{\nabh}$, $a,b \in \Nte$, where we use the
short-hand notation $\cN_{\nabh}:= \cN_{\Tr_{\nabh}}$ and where
$\Tr_{\nabh}$ is the \nsf weight on $B(H)$ with density $\nabh$. Also,
the elements $(y \ot 1) \al(b)$ span a core for the GNS-map of
$\tetiltil$.

On $B(H) \ot N$, we have two natural weights: $\tetiltil$ and $\Tr \ot
\te$. The following result was first proven by Yamanouchi in
\cite{Yam1} and called the Takesaki duality theorem for weights. We
indicate how it follows from \cite{V}.

\begin{lemma}
The following formula holds:
$$[D \tetiltil : D (\Tr_{\nabh} \ot \te)]_t = \nabtil_\te^{it} (\nabh^{-it}
\ot \nabte^{-it}) \; .$$
\end{lemma}
\begin{proof}
Writing $T_{\alhhklein}$ for the operator valued weight from $B(H) \ot N$
to $M \kruisje{\al} N$ associated with the integrable action $\alhh$,
we know from \cite{V}, Proposition 5.7 that $T_{\alhhklein}(z) = \Jtil_\te
T_{\alh}^{-1}(\Jtil_\te z \Jtil_\te) \Jtil_\te$. Here, $T_{\alh}$ is the operator
valued weight from $M \kruisje{\al} N$ to $\al(N)$ and $T_{\alh}^{-1}$
is the commutant operator valued weight from $\al(N)'$ to $(M \kruisje{\al} N)'$.
Using the language of spatial
derivatives, this means that
$$\Bigl( \frac{d \tetiltil}{d(1 \ot \te(\Jte \cdot \Jte))} \Bigr)^{it}
= \Bigl( \frac{d (\tetil(\Jtil_\te \cdot \Jtil_\te) \na T_{\alh}^{-1})}{d(\te
  \na \al^{-1})} \Bigr)^{-it} =
\Bigl( \frac{d (\tetil(\Jtil_\te \cdot \Jtil_\te))}{d(\te
  \na \al^{-1} \na T_{\alh})} \Bigr)^{-it} = \Bigr( \frac{d (\tetil(\Jtil_\te
  \cdot \Jtil_\te))}{d \tetil} \Bigr)^{-it} = \nabtil_\te^{it} \; .$$
Then, we get that
$$[D \tetiltil : D (\Tr_{\nabh} \ot \te)]_t = \Bigl( \frac{d
  \tetiltil}{d(1 \ot \te(\Jte \cdot \Jte))} \Bigr)^{it}
\Bigl( \frac{d (\Tr_{\nabh} \ot \te)}{d(1 \ot \te(\Jte \cdot \Jte))} \Bigr)^{-it}
= \nabtil_\te^{it} (\nabh^{-it} \ot \nabte^{-it}) \; .$$
This concludes the proof.
\end{proof}

\begin{definition} \label{def.radnik}
Define $[D \te \na \al : D \te]_t := \nabtil_\te^{it} (\nabh^{-it} \ot
\nabte^{-it})$. We call $[D \te \na \al : D \te]_t$ the Radon-Nikodym
derivative of $\te$ under the action $\al$ of $(M,\de)$.
\end{definition}

We now prove the main properties of $[D \te \na \al : D \te]_t$.

\begin{theorem} \label{theorem.cocycle}
Write $D_t := [D \te \na \al : D \te]_t$.
\begin{itemize}
\item $D_t \in M \ot N$.
\item $(\de \ot \io)(D_t) = (\io \ot \al)(D_t) (1 \ot D_t)$.
\end{itemize}
\end{theorem}
\begin{proof}
Denote by $(\sitiltil)_t$ the modular automorphism group of the bidual
weight $\tetiltil$ on $B(H) \ot N$. Because the dual weight $\tetil$
is $\sdeh$-invariant (see \cite{V}, Proposition 2.5), we conclude from the
proof of \cite{V}, Proposition 4.3, that $\sitiltil_t(x \ot 1) =
\sdeh^{-it} \nab^{-it} x \nab^{it} \sdeh^{it} \ot 1$ for all $x \in
M'$. From \cite{VVD}, Proposition 2.4, we know that $\sdeh^{-it}
\nab^{-it} = \nabh^{it} \sde^{it}$ and hence, $\sitiltil_t(x \ot 1) =
\nabh^{it} x \nabh^{-it} \ot 1$ for all $x \in
M'$. We easily conclude that $D_t \in M \ot N$.

Next, we introduce the notation $\te^\al := \tetiltil$ to stress that
we take a bidual weight with respect to the action $\al$. Consider the
action $\be:= (\si \ot \io)(\io \ot \al)$ of $(M,\de)$ on $B(H) \ot
N$. Because $\te^\al$ is an \nsf weight on $B(H) \ot N$, we can define
its bidual weight $(\te^\al)^\be$ on $B(H \ot H) \ot N$. From the
discussion above, we know a canonical GNS-map for $(\te^\al)^\be$
given by $$(y \ot z \ot 1) \be(\al(x)) \mapsto \la_{\nabh}(y) \ot
\la_{\nabh}(z) \ot \late(x) \quad\text{for all}\quad y,z \in
\cN_{\nabh}, x \in \Nte \; ,$$ where $\la_{\nabh}$ is a GNS-map for
$\Tr_{\nabh}$. We also know that the elements $(y \ot z \ot 1)
\be(\al(x))$ span a core for the corresponding GNS-map. Observe that
$(\Si V^* \Si) (\nabh \ot \nabh) (\Si V \Si) = Q \ot \nabh$, where $Q$
is the closure of $\sde \nabh$. Also observe that $\be(\al(x)) =
(\deop \ot \io)\al(x)$ for all $x \in N$. Using the above core for
$(\te^\al)^\be$, we conclude that
$$(\te^\al)^\be \na \Ad(\Si V \Si \ot 1) = \Tr_Q \ot \te^\al
\quad\text{and}\quad (\Tr_{\nabh} \ot \Tr_{\nabh} \ot \te) \na \Ad(\Si
V \Si \ot 1) = \Tr_Q \ot \Tr_{\nabh} \ot \te \; .$$
So, on the one hand, we have,
$$
[D ((\te^\al)^\be) : D (\Tr_{\nabh} \ot \Tr_{\nabh} \ot \te)]_t =
(\Si V \Si \ot 1) [D (\Tr_Q \ot \te^\al) : D (\Tr_Q \ot \Tr_{\nabh}
\ot \te)]_t (\Si V^* \Si \ot 1) = (\deop \ot \io)(D_t) \; .$$
On the other hand, we have
\begin{align}
[D ((\te^\al)^\be) & : D (\Tr_{\nabh} \ot \Tr_{\nabh} \ot \te)]_t =
[D ((\te^\al)^\be) : D( (\Tr_{\nabh} \ot \te)^\be) ]_t \; [D(
(\Tr_{\nabh} \ot \te)^\be) : D (\Tr_{\nabh} \ot \Tr_{\nabh} \ot
\te)]_t \notag \\ &= \be\bigl( [D \te^\al : D (\Tr_{\nabh} \ot \te)]_t \bigr) \;
[D ( (\Tr_{\nabh} \ot \te^\al)(\si \ot \io)) : D (\Tr_{\nabh} \ot
\Tr_{\nabh} \ot \te)]_t \label{calculation} \\ &= (\si \ot \io) \bigl( (\io \ot \al)(D_t) (1
\ot D_t) \bigr) \; , \notag
\end{align}
where we used the following well known fact from the theory of
operator valued weights: $[D \tetil : D \tilde{\mu}]_t = \be([D \te : D
\mu]_t)$ where $\be$ is an action of a l.c.\ quantum group and
$\te,\mu$ are \nsf weights with dual weights $\tetil,\tilde{\mu}$.
Combining both calculations above, the proof is finished.
\end{proof}

The proof of the previous theorem can be adapted to prove the
following useful and seemingly obvious formula.
\begin{lemma} \label{lemma.formula}
Suppose that $N$ is a
von Neumann algebra with \nsf weight $\te$ and let $\al$ and $\be$ be a actions
of $(M_1,\deo)$, resp.\ $(M_2,\det)$ on $N$. Suppose that there exists
a $^*$-isomorphism $\tw : M_1 \ot M_2 \recht M_1 \ot M_2$ ruling, in a
sense, the commutation between $\al$ and $\be$. More precisely,
$$(\io \ot \be)\al = (\tw\si \ot \io) (\io \ot \al) \be \;
.$$
If $\tw (\tauo_t \ot \taut_t) = (\tauo_t \ot \taut_t) \tw$, then
$$(\io \ot \be)\bigl( [D \te \na \al : D \te ]_t \bigr) \; (1 \ot [D \te
\na \be : D \te]_t) = (\tw \si \ot \io) \Bigl( (\io \ot \al)\bigl([D
\te \na \be : D \te]_t \bigr) \; (1 \ot [D \te \na \al : D \te]_t) \Bigr)
\; .$$
\end{lemma}
\begin{proof}
Consider the amplified action $\be_a:=(\si \ot \io)(\io \ot \be)$ of
$(M_2,\det)$ on $B(H_1) \ot N$. First, we have the bidual weight
$\te^\al$ on $B(H_1) \ot N$, which has again a bidual weight
$(\te^\al)^{\be_a}$ on $B(H_2 \ot H_1) \ot N$.
The same calculation as in Equation~\eqref{calculation} yields that
\begin{equation} \label{eq.teal}
[D (\te^\al)^{\be_a} : D (\Tr_{\nabh_2} \ot \Tr_{\nabh_1} \ot
\te)]_t = (\si \ot \io) \Bigl( (\io \ot \be)\bigl( [D \te \na \al : D
\te ]_t \bigr) \; (1 \ot [D \te \na \be : D \te]_t) \Bigr) \; .
\end{equation}
By symmetry, we have an amplified action $\al_a$ and
\begin{equation} \label{eq.tebe}
[D (\te^\be)^{\al_a} : D (\Tr_{\nabh_1} \ot \Tr_{\nabh_2} \ot
\te)]_t = (\si \ot \io) \Bigl( (\io \ot \al)\bigl( [D \te \na \be : D
\te ]_t \bigr) \; (1 \ot [D \te \na \al : D \te]_t) \Bigr) \; .
\end{equation}
Consider the unbounded operator $K := J_1 \sde_1 J_1 \ot J_2 \sde_2
J_2$, affiliated with $M_1' \ot M_2'$. Then, $K \ot 1$ is invariant
under the modular automorphism group of $\Tr_{\nabh_1} \ot \Tr_{\nabh_2} \ot
\te$ and it commutes with $$[D (\te^\be)^{\al_a} : D (\Tr_{\nabh_1} \ot
\Tr_{\nabh_2} \ot \te)]_t \in M_1 \ot M_2 \ot N \; .$$ Hence, $K \ot 1$ is invariant under the
modular automorphism group of $(\te^\be)^{\al_a}$. So, we can define the
\nsf weight $\mu:=((\te^\be)^{\al_a})_{K \ot 1}$. We know a GNS-map
for $(\te^\be)^{\al_a}$ and because the
closure of $(\nabh_1 \ot \nabh_2)K$ is $P_1 \ot P_2$, we have a
canonical GNS-map for $\mu$ given by
$$(z \ot 1) \al_a(\be(x)) \mapsto \la_{P_1 \ot P_2}(z) \ot \late(x)
\quad\text{for all}\quad z \in \cN_{P_1 \ot P_2}, x \in \Nte \; .$$
The elements $(z \ot 1) \al_a(\be(x))$ span a core for this
GNS-map. Denote by $\cZ$ the canonical implementation of the
automorphism $\tw$, which is a unitary on $H_1 \ot H_2$. Because $\tw$
commutes with $\tauo_t \ot \taut_t$ and because $P_1^{it} \ot
P_2^{it}$ is the canonical implementation of $\tauo_t \ot \taut_t$, we
get that $\cZ$ commutes with $P_1 \ot P_2$. Using the core for the
GNS-map of $(\te^\be)^{\al_a}$ and using the given formula
$$(\io \ot \be)\al = (\tw\si \ot \io) (\io \ot \al) \be \; ,$$ we conclude
that
$$\mu = \rho \na \Ad(\Si \cZ \ot 1) \; ,$$
where $\rho:= ((\te^\al)^{\be_a})_{\Si K \Si \ot 1}$.

Finally, it follows from Equation~\eqref{eq.tebe} that
$$[D \mu : D (\Tr_{P_1 \ot P_2} \ot \te)]_t = (\si \ot \io) \Bigl( (\io \ot \al)\bigl( [D \te \na \be : D
\te ]_t \bigr) \; (1 \ot [D \te \na \al : D \te]_t) \Bigr) \; ,$$
while it follows from Equation~\eqref{eq.teal} that
$$[D \rho :   D (\Tr_{P_2 \ot P_1} \ot \te)]_t = (\si \ot \io) \Bigl( (\io \ot \be)\bigl( [D \te \na \al : D
\te ]_t \bigr) \; (1 \ot [D \te \na \be : D \te]_t) \Bigr) \; .$$
Applying $\Ad(\Si \cZ \ot 1)$ and using that $\mu = \rho \na \Ad(\Si
\cZ \ot 1)$ and that $\cZ$ commutes with $P_1 \ot P_2$, we arrive at
the statement of the lemma.
\end{proof}

As a corollary, we prove the following result on closed quantum
subgroups.

\begin{proposition}
Let $(M,\de)$ be a l.c.\ quantum group. Suppose that $N \subset M$ is
a von Neumann subalgebra such that $\de(N) \subset N \ot N$.

Then, $(N,\de|_N)$ is a l.c.\ quantum group if and only if $R(N)
= N$ and $\tau_t(N) = N$ for all $t \in \R$.
\end{proposition}
\begin{proof}
If $(N,\de|_N)$ is a l.c.\ quantum group, the inclusion map is a
morphism and we know that $R(N)
= N$ and $\tau_t(N) = N$ for all $t \in \R$. So, we only have to prove the converse implication.

Define $\al:=\de|_N$, which we consider as a left action of $(M,\de)$
on $N$ and $\be:=\deop|_N$, which we consider as a left action of
$(M,\deop)$ on $N$. Choose an n.s.f.\ weight $\te$ on $N$ and define
$$u_t := [D \te \circ \al : D \te]_t \in M \ot N \; , \quad w_t := [D
\te \circ \be : D \te]_t \in M \ot N \; .$$
We will show that there exists a $(\si^\te_t)$-cocycle $(v_t)$ in $N$
such that $u_t = \al(v_t^*)(1 \ot v_t)$ for all $t \in \R$. If we then
correct the weight $\te$ with the cocycle $(v_t)$, we find a left
invariant weight on $(N,\de|_N)$.

Because $\al(N) \subset N \ot N$, we can use the argument of the proof
of Proposition 3.12 in \cite{V} to conclude that $u_t \in N \ot N$ and
analogously, $w_t \in N \ot N$ for all $t \in \R$. Remark that
at this point, we use in a crucial way that $R$ and $\tau_t$ leave $N$
globally invariant.

By Theorem
\ref{theorem.cocycle}, we know that $(\de \ot \io)(u_t) = (\io \ot
\de)(u_t) (1 \ot u_t)$, which we can rewrite as
$$(\de \ot \io)(\util_t W) = (\util_t W)_{13} (\util_t W)_{23}
\quad\text{with}\quad \util_t = (\Jh \ot J)u_t^*(\Jh \ot J) \; .$$
Because $\util_t \in M \ot M'$, the corepresentation $\util_t W$ is
covariant with respect to the standard representation of $M$ on $H$ in
the sense that $\de(x) = (\util_t W)^* (1 \ot x) (\util_t W)$ for all
$x \in M$. As we remarked after Definition 2.4 in \cite{BSV}, this
gives the existence of faithful, normal $^*$-homomorphisms $\pih_t :
\Mh \recht \B(H)$ satisfying $(\io \ot \pih_t)(W) = \util_t W$.

On the other hand, we have analogously,
$$(\deop \ot \io)(\wtil_t \; \Si V^* \Si) = (\wtil_t  \; \Si V^* \Si)_{13}
(\wtil_t  \; \Si V^* \Si)_{23} \quad\text{where}\quad \wtil_t = (\Jh \ot
J) w_t^* (\Jh \ot J) \; .$$
So, $\wtil_t  \; \Si V^* \Si$ is a corepresentation of $(M,\deop)$ which
is covariant with respect to the standard representation of
$(M,\deop)\hat{\text{ }} = \Mh'$. For the same reason as above, we
find faithful, normal $^*$-homomorphisms $\pih_t' :
\Mh' \recht \B(H)$ satisfying $(\io \ot \pih_t')(\Si V^* \Si) =
\wtil_t  \; \Si V^* \Si$.

Applying Lemma \ref{lemma.formula} to the actions $\al$ and $\be$
(with $\tw = \io$), we find that
$$(\io \ot \de)(w_t) (1 \ot u_t) = (\si \ot \io)\bigl( (\io \ot
\deop)(u_t) (1 \ot w_t) \bigr) \; ,$$
which can be rewritten as $$(\Si V^* \Si \;  w_t)_{13}
\quad\text{and}\quad (W u_t)_{23} \quad\text{commute}\; .$$
This means that the ranges of $\pih_t$ and $\pih_t'$ commute.

Write $\pi(x) = x$ for $x \in M$. Define $\cU_t := (\pih_t' \ot
\io)(V^*) \; (\pi \ot \io)(W^*) \; (\pih_t \ot \io)(\Vtil)$ with $\Vtil
= (J \ot J)\Wh^*(J \ot J)$. Checking it on a slice of $W$, we get that
$(\pih_t \ot \io)\dehop(x) = (\pi \ot \io)(W) (\pih_t(x) \ot 1) (\pi
\ot \io)(W^*)$ for all $x \in \Mh$. So,
$$(\pi \ot \io)(W^*) \; (\pih_t \ot \io)(\Vtil) \; (1 \ot a) =
(\pih_t(a) \ot 1) (\pi \ot \io)(W^*) \; (\pih_t \ot \io)(\Vtil)
\quad\text{for all}\quad a \in \Mh \; .$$
So, we arrive at
$$\cU_t (1 \ot a) \cU_t^* = \pih_t(a) \ot 1 \quad\text{for}\quad a \in
\Mh \; , \qquad \cU_t(1 \ot b) \cU_t^* = \pi(b) \ot 1
\quad\text{for}\quad b \in M \; ,$$
where the second formula is checked in an analogous way as the first one.
Because $\Mh M$ is strongly dense in $\B(H)$, it follows that there
exist normal, faithful, $^*$-homomorphisms $\Psi_t : \B(H) \recht
\B(H)$ such that $\cU_t(1 \ot x) \cU_t^* = \Psi_t(x) \ot 1$ for all $x
\in \B(H)$. It is also clear that the map $t \mapsto \Psi_t(x)$ is
strong$^*$-continuous for all $x \in \B(H)$.

We prove that $\Psi_t(\B(H)) = \B(H)$. So, suppose that $x \in
\Psi_t(\B(H))'$. Because $\Psi_t(a) = a$ for $a \in M$, it follows
that $x \in M'$ and we write $y = JxJ \in M$. Further, $x \in
\pih_t(\Mh)'$, which means that $1 \ot x$ commutes with $\util_t
W$. Hence, $1 \ot y$ commutes with $W u_t$. So,
$$\de(y) = u_t(1 \ot y)u_t^* \; .$$
From this, we conclude that $\deop(y) \in M \ot N$ and obviously
$(\deop \ot \io)\deop(y) = (\io \ot \be)\deop(y)$. This gives us that
$\deop(y) \in \be(N)$ and hence, $y \in N$. So,
$$\al(y) = u_t(1 \ot y) u_t^* = \nabtil^{it} (1 \ot \si^\te_{-t}(y))
\nabtil^{-it} \; ,$$ from which we conclude that $\de(\si^\te_{-t}(y))
= 1 \ot \si^\te_{-t}(y)$. We finally find that $\si^\te_{-t}(y) \in
\C$ and so, $y \in \C$. This proves that $\Psi_t(\B(H)) = \B(H)$.

So, $\Psi_t$ is a family of isomorphisms of $\B(H)$ which is pointwise
strong$^*$-continuous. It is then easy to find a strongly continuous family $(\vtil_t)$
of unitaries in $\B(H)$ such that $\Psi_t = \Ad \vtil_t$ for all $t
\in \R$.
Because $\Psi_t(a) = a$ for all $a \in M$, we get $\vtil_t \in
M'$. We define $v_t := J \vtil_t^* J$ and conclude from the equation
$(\io \ot \pih_t)(W) = \util_t W$ that
$$u_t^* W^* = (1 \ot v_t^*)W^* (1 \ot v_t) \quad\text{and hence,}\quad
u_t = \de(v_t^*)(1 \ot v_t) \; .$$
But then, we observe that $\de(v_t) \in N \ot M$, because $u_t \in N
\ot N$, and the same reasoning as above yields $v_t \in N$ for all $t
\in \R$. So, we write $u_t = \al(v_t^*) (1 \ot v_t)$.

Because $u_t = [D \tetiltil : D (\Tr_{\nabh} \ot \te)]_t$, we know
that $u_{t+s} = u_t \; (\tau_t \ot \si^\te_t)(u_s)$. So, we can
calculate
\begin{align*}
\al(v_{t+s}^*)(1 \ot v_{t+s}) & = u_{t+s} = u_t \; (\tau_t \ot
\si^\te_t)(\al(v_s^*)(1 \ot v_s)) = \nabtil^{it} (\nabh^{-it} \ot
\nabte^{-it}) \; (\tau_t \ot
\si^\te_t)(\al(v_s^*)(1 \ot v_s)) \\ &= \nabtil^{it} \al(v_s^*)(1\ot
v_s) (\nabh^{-it} \ot
\nabte^{-it}) = \al(\si^\te_t(v_s^*)) \; u_t \; (1 \ot \si^\te_t(v_s))
\\ &= \al( \si^\te_t(v_s^*) v_t^*) (1 \ot v_t \si^\te_t(v_s)) \; .
\end{align*}
It follows that $v_t \si^\te_t(v_s) v_{t+s}^*$ is invariant under
$\de$ and hence, a scalar. We find a continuous function $\sla : \R^2
\recht \T$ such that $v_{t+s} = \sla(t,s) v_t \si^\te_t(v_s)$. It is
then easy to check that $\sla$ is a $2$-cocycle and so, necessarily, a
$2$-coboundary. This means that there exists a continuous function
$\mu : \R \recht \T$ such that $\sla(t,s) = \mu(t) \mu(s)
\overline{\mu(t+s)}$. If we replace $v_t$ by $\mu(t) v_t$, we still
have $u_t = \al(v_t^*)(1 \ot v_t)$ and we get moreover that $v_{t+s} =
v_t \si^\te_t(v_s)$ for all $s,t \in \R$. Because the map $t \mapsto
v_t$ is strong$^*$-continuous, it follows that $(v_t)$ is a
$(\si^\te_t)$-cocycle.

So, there exists a unique n.s.f.\ weight $\vfi_N$ on $N$ such that $[D
\vfi_N : D \te]_t = v_t$. Then,
\begin{align*}
[D \vfitiltil_N : D (\Tr_{\nabh} \ot \vfi_N)]_t &=
[D \vfitiltil_N : D \tetiltil]_t \; [D \tetiltil : D (\Tr_{\nabh} \ot
\te)]_t \; [D (\Tr_{\nabh} \ot \te) : D (\Tr_{\nabh} \ot \vfi_N)]_t \\
&= \al(v_t) \; u_t \; (1 \ot v_t^*) = 1 \; .
\end{align*}
So, $[D \vfi_N \circ \al : D \vfi_N]_t = 1$, which means that $(\io
\ot \vfi_N)\al(x) = \vfi_N(x) \; 1$ for all $x \in
\cM_{\vfi_N}^+$. Hence, $\vfi_N$ is a left invariant weight on
$(N,\de|_N)$.

Because $R(N) = N$, we can restrict $R$ to $N$ and we find an
anti-automorphism of $N$ that anti-commutes with $\de|_N$. So, we also
get a right invariant weight and $(N,\de|_N)$ is a l.c.\ quantum group.
\end{proof}


\begin{thebibliography}{AA}

\bibitem{B} {\sc S. Baaj}, Repr{\'e}sentation r{\'e}guli{\`e}re du groupe quantique des d{\'e}placements de Woronowicz.
{\it Ast{\'e}risque} {\bf 232} (1995), 11--48.

\bibitem{BS}  {\sc S. Baaj  \&  G. Skandalis},
Unitaires multiplicatifs et dualit{\'e} pour les produits crois{\'e}s de C$^*$-alg{\`e}bres. {\it
Ann. Scient. Ec. Norm. Sup., $4{}^e$ s{\'e}rie}, {\bf 26} (1993), 425--488.

\bibitem{BSV} {\sc S. Baaj, G. Skandalis \& S. Vaes}, Non-semi-regular
  quantum groups coming from number theory. {\it Comm. Math. Phys.},
  to appear.

\bibitem{Drin} {\sc V. G. Drinfel'd}, Quantum groups. {\it In}
Proceedings ICM (Berkeley, Calif., 1986), Amer. Math. Soc.,
Providence, RI, 1987, pp. 798--820.

\bibitem{E-S} {\sc M. Enock \& J-M. Schwartz}, Kac Algebras and Duality of Locally Compact Groups, Springer-Verlag, 1992.

\bibitem{haagerup} {\sc U. Haagerup}, Operator valued weights in von
Neumann algebras I. {\it J. Funct. Anal.} {\bf 32} (1979),
175--206.

\bibitem{JK} {\sc J. Kustermans}, Locally compact quantum groups in the universal setting.
{\it Int. J. Math.} {\bf 12} (2001), 289--338.

\bibitem{KV1} {\sc J. Kustermans \& S. Vaes}, Locally compact
quantum groups. {\it Ann. Scient. Ec. Norm. Sup., $4{}^e$ s{\'e}rie} {\bf 33} (2000), 837--934.

\bibitem{KV2} {\sc J. Kustermans \& S. Vaes}, Locally compact quantum groups in the von Neumann
algebraic setting. {\it Math. Scand.}, to appear.

\bibitem{Maj} {\sc S. Majid}, More examples of bicrossproduct and double cross product Hopf
algebras. {\it Isr. J. Math.} {\bf 72} (1990), 133--148.

\bibitem{Mas-Nak} {\sc T. Masuda \& Y. Nakagami}, A von Neumann algebraic framework for the duality of the quantum groups. {\it Publ. RIMS, Kyoto University} {\bf 30} (1994), 799--850.

\bibitem{V3} {\sc S. Vaes}, Strictly outer actions of locally compact groups and
  quantum groups. {\it In preparation.}

\bibitem{V} {\sc S. Vaes}, The unitary implementation of a
locally compact quantum group action. {\it J. Func. Anal.} {\bf
180} (2001), 426--480.

\bibitem{V2} {\sc S. Vaes}, A Radon-Nikodym theorem for von Neumann algebras.
       {\it J. Operator Theory} {\bf 46} (3) (2001), 477--489.

\bibitem{VV} {\sc S. Vaes \& L. Vainerman}, Extensions of locally
compact quantum groups and the bicrossed product construction. {\it
  Adv. in Math.}, to appear.

\bibitem{VVD} {\sc S. Vaes \& A. Van Daele}, The Heisenberg
  commutation relations, commuting squares and the Haar measure on
  locally compact quantum groups. {\it Proceedings of the OAMP
  Conference, Constantza, 2001}, to appear.

\bibitem{VD} {\sc A. Van Daele}, An algebraic framework for group duality. {\it Adv. in Math.} {\bf 140} (1998), 323--366.

\bibitem{Wor1} {\sc S.L. Woronowicz}, Compact matrix
  pseudogroups. {\it Comm. Math. Phys.} {\bf 111} (4) (1987), 613--665.

\bibitem{Wor2} {\sc S.L. Woronowicz}, Compact quantum groups. {\it In}
  Sym{\'e}tries quantiques (Les Houches, 1995), North-Holland, Amsterdam,
  1998, pp. 845--884.

\bibitem{Wor3} {\sc  S.L. Woronowicz}, From multiplicative unitaries to
  quantum groups. {\it Int. J. Math.} {\bf 7} (1) (1996), 127--149.

\bibitem{Yam1} {\sc T. Yamanouchi}, Takesaki duality for weights on
  locally compact quantum group covariant systems. {\it J. Operator
  Theory}, to appear.

\bibitem{Yam2} {\sc T. Yamanouchi}, Canonical extension of actions of
  locally compact quantum groups. {\it Preprint}.

\bibitem{Yam3} {\sc T. Yamanouchi},
Double group construction of quantum groups in the von Neumann algebra
framework.  {\it J. Math. Soc. Japan} {\bf 52} (4) (2000), 807--834.

\end{thebibliography}
\end{document}